\documentclass[12pt]{article}
\usepackage{amsmath,amsfonts,amssymb,amsthm}
\input amssym.def
\begin{document}
\def \Z{\mathbb Z}
\def \C{\mathbb C}
\def \R{\Bbb R}
\def \Q{\mathbb Q}
\def \N{\mathbb N}

\def \A{{\mathcal{A}}}
\def \D{{\mathcal{D}}}
\def \E{{\mathcal{E}}}
\def \E{{\mathcal{E}}}
\def \H{\mathcal{H}}
\def \S{{\mathcal{S}}}
\def \wt{{\rm wt}}
\def \tr{{\rm tr}}
\def \span{{\rm span}}
\def \Res{{\rm Res}}
\def \Der{{\rm Der}}
\def \End{{\rm End}}
\def \Ind {{\rm Ind}}
\def \Irr {{\rm Irr}}
\def \Aut{{\rm Aut}}
\def \GL{{\rm GL}}
\def \Hom{{\rm Hom}}
\def \mod{{\rm mod}}
\def \ann{{\rm Ann}}
\def \ad{{\rm ad}}
\def \rank{{\rm rank}\;}
\def \<{\langle}
\def \>{\rangle}

\def \g{{\frak{g}}}
\def \h{{\hbar}}
\def \k{{\frak{k}}}
\def \sl{{\frak{sl}}}
\def \gl{{\frak{gl}}}

\def \be{\begin{equation}\label}
\def \ee{\end{equation}}
\def \bex{\begin{example}\label}
\def \eex{\end{example}}
\def \bl{\begin{lem}\label}
\def \el{\end{lem}}
\def \bt{\begin{thm}\label}
\def \et{\end{thm}}
\def \bp{\begin{prop}\label}
\def \ep{\end{prop}}
\def \br{\begin{rem}\label}
\def \er{\end{rem}}
\def \bc{\begin{coro}\label}
\def \ec{\end{coro}}
\def \bd{\begin{de}\label}
\def \ed{\end{de}}

\newcommand{\m}{\bf m}
\newcommand{\n}{\bf n}
\newcommand{\nno}{\nonumber}
\newcommand{\nord}{\mbox{\scriptsize ${\circ\atop\circ}$}}
\newtheorem{thm}{Theorem}[section]
\newtheorem{prop}[thm]{Proposition}
\newtheorem{coro}[thm]{Corollary}
\newtheorem{conj}[thm]{Conjecture}
\newtheorem{example}[thm]{Example}
\newtheorem{lem}[thm]{Lemma}
\newtheorem{rem}[thm]{Remark}
\newtheorem{de}[thm]{Definition}
\newtheorem{hy}[thm]{Hypothesis}
\makeatletter \@addtoreset{equation}{section}
\def\theequation{\thesection.\arabic{equation}}
\makeatother \makeatletter

\begin{center}
{\Large \bf Regular representations and $A_{n}(V)$-$A_{m}(V)$ bimodules}
\end{center}

\begin{center}
{Haisheng Li\\
Department of Mathematical Sciences\\
Rutgers University, Camden, NJ 08102}
\end{center}

\begin{abstract}
This paper is to establish a natural connection between regular representations for a vertex operator algebra $V$ and
$A_{n}(V)$-$A_{m}(V)$ bimodules of Dong and Jiang. Let $W$ be a weak $V$-module and let  
$(m,n)$ be a pair of nonnegative integers. We study 
two quotient spaces $A_{n,m}^{\dagger}(W)$ and $A^{\diamond}_{n,m}(W)$ of $W$. 
 It is proved that the dual space $A^{\dagger}_{n,m}(W)^{*}$ 
viewed as a subspace of $W^*$ coincides with the level-$(m,n)$ vacuum subspace 
of the regular representation module $\mathfrak{D}_{(-1)}(W)$. 
By making use of this connection, we obtain an $A_{n}(V)$-$A_m(V)$ bimodule structure 
on both $A_{n,m}^{\dagger}(W)$ and $A^{\diamond}_{n,m}(W)$. Furthermore, we obtain an $\N$-graded weak $V$-module structure 
together with a commuting right $A_m(V)$-module structure on $A^{\diamond}_{\Box,m}(W):=\oplus_{n\in \N}A^{\diamond}_{n,m}(W)$. 
Consequently, we recover the corresponding results and roughly confirm a conjecture  of Dong and Jiang.
\end{abstract}

\section{Introduction}
For a vertex operator algebra $V$, what was called the regular representation (module),
associated to a weak $V$-module $W$ and a nonzero complex number $z$, 
is a weak $V\otimes V$-module $\mathfrak{D}_{(z)}(W)$
which was constructed canonically inside (the full dual space) $W^*$ (see \cite{Li-reg}). 
In case $W=V$ (the adjoint module), it was proved that 
the socle of the weak $V\otimes V$-module $\mathfrak{D}_{(z)}(V)$ admits a Peter-Weyl type decomposition.

It turned out that regular representation has deep connections with various theories in the field of vertex operator algebras.
First of all, it has an intrinsic connection with Huang and Lepowsky's tensor product theory (see \cite{hl-1, hl-3}, \cite{huang-4}).
In this very theory, for any $V$-modules $W_1$, $W_2$ and for any nonzero complex number $z$,
a $V$-module is constructed inside the full dual space $(W_1\otimes W_2)^{*}$ and its contragredient module
 is defined to be the tensor product of $W_1$ and $W_2$. 
 A result of \cite{Li-reg} states that for $V$-modules $W_1, W_2$ and $W$, 
a linear map $F: W_1\otimes W_2\rightarrow W^{*}$ is a
$P(z)$-intertwining map in the sense of Huang and Lepowsky if and only if $F(W_1\otimes W_2)\subset \mathfrak{D}_{(z)}(W)$ 
and $F$ is a $V\otimes V$-module homomorphism from $W_1\otimes W_2$ to $\mathfrak{D}_{(z)}(W)$.
(It was proved in \cite{hl-2} that a $P(z)$-intertwining map amounts to an intertwining operator in the sense of \cite{fhl}.)
The relationship between regular representation and Huang-Lepowsky tensor functor 
was explored further in \cite{Li-reg-hl}.

Regular representation also has natural connections with Zhu's $A(V)$-theory (see \cite{zhu1, zhu2}) 
and its generalization $A_n(V)$-theory (see \cite{dlm-anv}). 
The essence of Zhu's $A(V)$-theory is that an associative algebra $A(V)$ is associated to each vertex operator algebra $V$ and 
 a natural bijection is established between the set of equivalence classes of irreducible $\N$-graded, namely admissible, 
 $V$-modules and the set of equivalence classes of irreducible $A(V)$-modules. 
Also in this theory, Frenkel-Zhu's fusion rule theorem gives a way to determine fusion rules 
by using $A(V)$-bimodules and $A(V)$-module homomorphisms (see \cite{fz}; cf. \cite{li-fusion}).
Zhu's $A(V)$-theory was generalized  in \cite{dlm-anv}, where
a sequence of associative algebras $A_n(V)$ for $n\ge 0$ was introduced with $A(V)=A_0(V)$.
This sequence of associative algebras is an important tool in the study on vertex operator algebra representations, 
and one of the results therein is that a vertex operator algebra 
$V$ is rational if and only if $A_n(V)$ for all $n\ge 0$ are (finite-dimensional) semisimple.

Connections of regular representation with Zhu's $A(V)$-theory and $A_n(V)$-theory 
were studied in \cite{Li-reg-fusion, Li-Anv}, where new proofs  for several known theorems were obtained. 
One of the main results states that for any $V$-module $W$, the dual of the Frenkel-Zhu $A(V)$-bimodule $A(W)$,
as a subspace of $W^{*}$, coincides with the vacuum space of $\mathfrak{D}_{(-1)}(W)$. 
This gives a direct connection of $A(W)$ with (weak) $V$-modules and intertwining operators, which leads to 
more conceptual and (arguably) easier proofs.
(See for example \cite{Li-reg-fusion}  for a different proof of the Frenkel-Zhu's fusion rule theorem; cf. \cite{li-fusion}.)

In this current work, our main concerns are $A_{n}(V)$-$A_{m}(V)$ bimodules which were introduced by Dong and Jiang.
Let $V$ be a vertex operator algebra. For $m,n\in \N$, 
Dong and Jiang in \cite{dj1} defined a quotient space $A_{n,m}(V)$ of $V$ and 
proved that $A_{n,m}(V)$ has a natural (left-right) $A_{n}(V)$-$A_{m}(V)$ bimodule structure. 
Among the main results, it was proved that for any $A_m(V)$-module $U$ with $m$ fixed, 
$\bigoplus_{n\in \N}A_{n,m}(V)\otimes_{A_m(V)}U$ 
has a canonical $V$-module structure with $A_{m,m}(V)\otimes_{A_m(V)}U=U$, which satisfies a certain universal property.
This theory generalizes the $A_n(V)$-theory in a natural way with new perspectives. 

The main goal of this paper is to establish a natural connection between regular representations and Dong-Jiang's 
$A_{n}(V)$-$A_{m}(V)$-bimodules for a vertex operator algebra $V$. 
For any weak $V$-module $W$ and for any pair $(m,n)$ of nonnegative integers, we study 
two quotient spaces  of $W$, denoted by $A_{n,m}^{\dagger}(W)$ and $A^{\diamond}_{n,m}(W)$. 
 It is proved that the dual space $A^{\dagger}_{n,m}(W)^{*}$ 
viewed as a subspace of $W^*$ coincides with the level-$(m,n)$ vacuum subspace 
of the regular representation module $\mathfrak{D}_{(-1)}(W)$. 
By making use of this connection, we obtain an $A_{n}(V)$-$A_m(V)$ bimodule structure 
on both $A_{n,m}^{\dagger}(W)$ and $A^{\diamond}_{n,m}(W)$. 
Furthermore, for any fixed $m\in \N$ we obtain an $\N$-graded weak $V$-module structure 
together with a commuting right $A_m(V)$-module structure on $A^{\diamond}_{\Box,m}(W):=\oplus_{n\in \N}A^{\diamond}_{n,m}(W)$. 
Using this we recover the aforementioned main result of \cite{dj1}.

We continue to describe the contents of this paper with some details.
Let $V$ be a vertex operator algebra and let $n\in \N$. From \cite{dlm-anv} 
we have an associative algebra $A_n(V)$ with an order $2$ anti-automorphism $\theta$.
 For any weak $V$-module $W$, set
$$\Omega_n(W)=\{ w\in W\ |\ x^{n}Y(x^{L(0)}v,x)w\in W[[x]]\ \ \text{ for }v\in V\}.$$
Then $\Omega_n(W)$ is naturally an $A_n(V)$-module.
Now, let $M$ be a weak $V\otimes V$-module. On $M$, we have two commuting $V$-module structures, 
denoted by $Y^{L}(\cdot,x)$ and $Y^R(\cdot,x)$, through the natural identifications 
$V=V\otimes \C{\bf 1}$ and $V=\C{\bf 1}\otimes V$.
For a pair $(m,n)$ of nonnegative integers, define the level-$(m,n)$ vacuum space
$$\Omega_{m,n}(M)=\Omega_{m}(M,Y^L)\cap \Omega_{n}(M,Y^R), $$
which is naturally an $A_m(V)\otimes A_n(V)$-module. 

Let $W$ be a weak $V$-module and let $m,n\in \N$. 
As the main objects of this paper, we study two quotient spaces 
$A_{n,m}^{\dagger}(W):=W/O^{\dagger}_{n,m}(W)$ and $A_{n,m}^{\diamond}(W):=W/O'_{n,m}(W)$, 
where $O^{\dagger}_{n,m}(W)$ is the subspace of $W$, linearly spanned by vectors
\begin{eqnarray*}
v\circ_{m}^{n}w:=\Res_{z}\frac{(1+z)^{m}}{z^{m+n+2}}Y((1+z)^{L(0)}v,z)w
\end{eqnarray*}
for $v\in V,\ w\in W$, and where
\begin{eqnarray*}
O'_{n,m}(W)=O^{\dagger}_{n,m}(W)+(L(-1)+L(0)+m-n)W,
\end{eqnarray*}
which is the same as in \cite{dj1}. As the first key result, we show that as subspaces of $W^{*}$,
\begin{eqnarray}
\Omega_{m,n}(\mathfrak{D}_{(-1)}(W))=(A^{\dagger}_{n,m}(W))^{*}.
\end{eqnarray}
View $(A^{\diamond}_{n,m}(W))^{*}$ naturally as  a subspace of $(A^{\dagger}_{n,m}(W))^{*}$ and define 
\begin{eqnarray}
\Omega^{\diamond}_{m,n}(\mathfrak{D}_{(-1)}(W))=(A^{\diamond}_{n,m}(W))^{*}.
\end{eqnarray}
It is proved that $\Omega^{\diamond}_{m,n}(\mathfrak{D}_{(-1)}(W))$ is an $A_m(V)\otimes A_n(V)$-submodule 
of $\Omega_{m,n}(\mathfrak{D}_{(-1)}(W))$. 
(More precisely, we use a deformed $A_m(V)\otimes A_n(V)$-module structure on $\Omega_{m,n}(\mathfrak{D}_{(-1)}(W))$, 
which corresponds to a particularly deformed $V\otimes V$ structure on $\mathfrak{D}_{(-1)}(W)$.)
Using this and the algebra anti-automorphism $\theta$, 
we then obtain an $A_n(V)$-$A_m(V)$-bimodule structure on
both $A^{\dagger}_{n,m}(W)$ and $A^{\diamond}_{n,m}(W)$.

 In \cite{dj1}, a subspace $O_{n,m}(V)$ of $V$ was introduced, which contains $O'_{n,m}(V)$,
 and it was proved that $A_{n,m}(V):=V/O_{n,m}(V)$ is an $A_n(V)$-$A_m(V)$-bimodule.
It was conjectured therein that $O_{n,m}(V)=O'_{n,m}(V)$ and this conjecture was confirmed in case $m=n$.
 While the actions of $A_n(V)$ in \cite{dj1} and in this current paper are the same, 
the right actions of $A_m(V)$ appear to be different.

Let $m\in \N$ be fixed. Consider the subspace $\Omega^{\diamond}_{m,\Box}(\mathfrak{D}_{(-1)}(W)):=
\sum_{n\in \N} \Omega^{\diamond}_{m,n}(\mathfrak{D}_{(-1)}(W))$ of $\mathfrak{D}_{(-1)}(W)$. It is shown that
\begin{eqnarray}
\Omega^{\diamond}_{m,\Box}(\mathfrak{D}_{(-1)}(W))=\bigoplus_{n\in \N} \Omega^{\diamond}_{m,n}(\mathfrak{D}_{(-1)}(W)),
\end{eqnarray} 
which is naturally an $\N$-graded weak $V$-module with a commuting (left) $A_m(V)$-module structure. 
On the other hand, follow \cite{dj1} to define an $\N$-graded vector space
\begin{eqnarray}
A^{\diamond}_{\Box,m}(W)=\bigoplus_{n\in \N}A^{\diamond}_{n,m}(W),
\end{eqnarray}
whose graded dual space coincides with $\Omega^{\diamond}_{m,\Box}(\mathfrak{D}_{(-1)}(W))$.
We then obtain  an $\N$-graded weak $V$-module structure together with a commuting
right $A_m(V)$-module structure on  $A^{\diamond}_{\Box,m}(W)$, 
with $\Omega^{\diamond}_{m,\Box}(\mathfrak{D}_{(-1)}(W))$ as its contragredient dual $V$-module.
This recovers a theorem of \cite{dj1}.  To a certain extent, this confirms Dong and Jiang's  conjecture.
The canonical connection established in this current work reconfirms from a different direction that 
$A_{n}(V)$-$A_{m}(V)$-bimodules are natural and important objects to study.

This paper is organized as follows: In Section 2, we review associative algebras $A_n(V)$, 
level-$n$ vacuum space $\Omega_n(W)$, and spaces $A^{\dagger}_{n,m}(W)$, $A^{\diamond}_{n,m}(W)$. 
 In Section 3, we recall some basic results on regular representations and 
 establish a canonical connection between $A_{n,m}^{\dagger}(W)$ and $\Omega_{m,n}(\mathfrak{D}_{(-1)}(W))$.
In Section 4, we study the connection between $\Omega^{\diamond}_{m,n}(\mathfrak{D}_{(-1)}(W))$ and $A^{\diamond}_{n,m}(W)$, and
study $\N$-graded $V$-modules $\Omega^{\diamond}_{m,\Box}(\mathfrak{D}_{(-1)}(W))$ and $A^{\diamond}_{\Box,m}(W)$.
 
 \vspace{0.5cm}
 
{\bf Acknowledgments:}  We thank Chongying Dong for extensive discussions.

\section{Associative algebra $A_n(V)$ and level-$n$ vacuum space $\Omega_n(W)$}

This section is preliminary;  We recall from \cite{dlm-anv} the sequence of associative algebras $A_n(V)$ and 
the level-$n$ vacuum space $\Omega_n(W)$ of a weak $V$-module $W$.
 Following Dong-Jiang \cite{dj1}, we introduce subspaces $O^{\dagger}_{n,m}(W)$  and $O'_{n,m}(W)$ of $W$ 
 for any pair $(m,n)$ of nonnegative integers and for any weak $V$-module $W$.

We begin by reviewing some basics on vertex operator algebras and modules (see \cite{flm}, \cite{fhl}). 
Let $V$ be a vertex operator algebra. 
The following are some of the main ingredients:  $V=\bigoplus_{n\in \Z}V_{(n)}$ is a $\Z$-graded 
vector space with a distinguished vector $\omega\in V_{(2)}$, called the {\em conformal vector}. 
Write $$Y(\omega,x)\ \left(=\sum_{n\in \Z}\omega_nx^{-n-1}\right)  =\sum_{n\in \Z}L(n)x^{-n-2}.$$
The following are the basic properties: 
$$[L(m),L(n)]=L(m+n)+\frac{1}{12}(m^3-m)\delta_{m+n,0}c$$
for $m,n\in \Z$, where $c\in \C$ is called the {\em central charge,} and  for $n\in \Z$,
\begin{eqnarray*}
&&V_{(n)}=\{ v\in V\ |\  L(0)v=nv\},\\
&&Y(L(-1)v,x)=\frac{d}{dx}Y(v,x)\quad \text{ for }v\in V.
\end{eqnarray*}
 A vector $v\in V_{(n)}$ with $n\in \Z$ is said to be {\em homogeneous}
 of {\em conformal weight $n$}  and we write $\wt v=n$. The {\em two grading restrictions} state
 \begin{eqnarray}
&& \dim V_{(n)}<\infty\quad \text{ for all }n\in \Z,\\
&&V_{(n)}=0\quad \text{ for all sufficiently negative integers }n.
 \end{eqnarray}
 
A $V$-module $W$ as a vector space is $\C$-graded with
$W=\bigoplus_{\lambda \in \C}W_{(\lambda)}$, where
$$W_{(\lambda)}=\{ w\in W\ |\ L(0)w=\lambda w\},$$
satisfying the two grading restrictions that  for $\lambda\in \C$,
$\dim W_{(\lambda)}<\infty$, and 
for any fixed $\lambda\in \C$, $W_{(\lambda+n)}=0$ for all sufficiently negative integers $n$.

By definition, a {\em weak $V$-module} is a module for $V$ viewed as a vertex algebra. 
Any weak $V$-module is naturally a module for the Virasoro algebra (of the same central charge), 
though the action of $L(0)$ might be nonsemisimple. 

The following are some basic facts from \cite{fhl} we shall need:
\begin{eqnarray}
&&x_0^{L(0)}Y(u,x)=Y(x_0^{L(0)}u,x_0x)x_0^{L(0)}, \label{L(0)-conjugation} \\
&&e^{-x_0L(1)}Y(u,x)= Y\left(e^{-x_0(1+x_0x)L(1)}(1+x_0x)^{-2L(0)}u, \frac{x}{1+x_0x}\right)e^{-x_0L(1)}. \quad
\label{L(1)-conjugation}
\end{eqnarray}
Also from \cite{fhl}, we have 
\begin{eqnarray}
x_1^{-L(0)}e^{xL(1)}x_1^{L(0)}&=&e^{xx_1L(1)},\\
(1+z_0x)^{-2L(0)}e^{z_0(1-z_0x)L(1)}&=&e^{z_0(1+z_0x)L(1)}(1-z_0x)^{2L(0)}.
\end{eqnarray}

Introduce a bijective linear operator on $V$
\begin{eqnarray}\label{def-theta}
\theta: \  V\rightarrow V; \  \  v\mapsto e^{L(1)}(-1)^{L(0)}v,
\end{eqnarray}
which plays an important role (see \cite{fhl}, \cite{zhu1}).
We have $\theta^2=1$ as
$$(-1)^{-L(0)}e^{L(1)}(-1)^{L(0)}=e^{-L(1)}\   \text{ and }\  (-1)^{-L(0)}=(-1)^{L(0)}\  \text{ on }V.$$

Next, we recall the sequence of associative algebras $A_n(V)$ from \cite{dlm-anv}.

\bd{def-O-n-W}
{\em Let $W$ be a weak $V$-module and let $n$ be a nonnegative integer.
For $v\in V, \ w\in W$, set
\begin{eqnarray}
v*_{n}w&=&\Res_{z}\sum_{i=0}^{n}\binom{-n-1}{i}\frac{(1+z)^{n}}{z^{n+1+i}}Y\left((1+z)^{L(0)}v,z\right)w,\label{*n-product}\\
v\circ_{n}w&=&\Res_{z}\frac{(1+z)^{n}}{z^{2n+2}}Y\left((1+z)^{L(0)}v,z\right)w.
\end{eqnarray}
Furthermore, define $O_{n}(W)$ to be the subspace of $W$, linearly spanned by vectors
$$(L(-1)+L(0))w,\   \   v\circ_{n}w \quad \text{(for all }v\in V,\ w\in W).$$}
\ed

The following is a result of \cite{dlm-anv}: 

\bt{AnV-algebra}
Let $n$ be any nonnegative integer. 
For the nonassociative algebra $(V,*_n)$ with the operation $*_n$ defined by (\ref{*n-product}), 
$O_n(V)$ is a two-sided ideal and the quotient algebra $V/O_n(V)$ is an associative algebra, denoted by 
$A_n(V)$. Furthermore, the linear operator $\theta=e^{L(1)}(-1)^{L(0)}$ on $V$
reduces to an anti-automorphism  of $A_n(V)$ with $\theta^2=1$. On the other hand,
the identity operator on $V$ gives rise to an algebra epimorphism $\psi_n: A_{n+1}(V)\rightarrow A_n(V)$ for every $n\in \N$.
\et

\bd{def-Omega-n}
{\em  Let $(W,Y_W)$ be a weak $V$-module, $n$ a nonnegative integer. Set
 \begin{eqnarray}
 \Omega_n(W)=\{ w\in W\ |\  x^nY_W(x^{L(0)}v,x)w\in W[[x]]\  \ \text{ for } v\in V  \}.
 \end{eqnarray}}
 \ed
 
 The following  was obtained in \cite{dlm-anv}, generalizing a result of Zhu (see \cite{zhu1}):
 
\bp{prop-Omega-n-W}
Let $(W,Y_W)$ be any weak $V$-module.  For $v\in V,\ w\in W$, define
 \begin{eqnarray}
 v\cdot w=\Res_x x^{-1}Y_W(x^{L(0)}v,x)w\in W.
 \end{eqnarray}
 Then 
$$u\cdot (v\cdot w)=(u*_nv)\cdot w \quad \text{for  }u,v\in V,\ w\in \Omega_n(W),$$ 
and $O_n(V)\cdot \Omega_n(W)=0$.
 Furthermore, $\Omega_n(W)$ is naturally an $A_n(V)$-module.
 \ep
 
 Recall the following technical result from
 \cite{Li-Anv} (Lemma 3.8\footnote{There is a typo for the two strict inequalities, which is fixed here}):
 
 \bl{lem3.8-anv}
 Let $W$ be a weak $V$-module, $n$ a nonnegative integer. Then for any (finitely many) homogeneous vectors 
 $v^{(1)},v^{(2)},\dots, v^{(r)}\in V$ and for any integers $m_i\in \Z$,
 \begin{eqnarray}
 v^{(1)}_{m_1}v^{(2)}_{m_2}\cdots v^{(r)}_{m_r}\cdot \Omega_n(W)=0
 \end{eqnarray}
 whenever $\wt \left(v^{(1)}_{m_1}v^{(2)}_{m_2}\cdots v^{(r)}_{m_r}\right)< -n$, i.e., 
 $$m_1+\cdots +m_r> \wt v^{(1)}+\cdots +\wt v^{(r)}-r+n.$$
 In particular, for any homogeneous vector $v\in V$,  $(v_m)^n\Omega_n(W)=0$ if $m\ge \wt v$.
 \el
 
 Note that for any weak $V$-module $W$, 
subspaces $\Omega_n(W)$ for $n\in \N$ form an ascending sequence.
 Using Lemma \ref{lem3.8-anv}, we immediately have:
 
 \bl{lemma-omega-union}
 Let $W$ be a weak $V$-module.  Set $\Omega_n(W)=0$ for $n<0$.  Then
   \begin{eqnarray}
  u_{m}\cdot \Omega_n(W)\subset \Omega_{n+\wt (u_m)} (W)
  \end{eqnarray}
for any homogeneous vector $u\in V$ and for any $m,n\in \Z$, where $\wt (u_m)=\wt u-m-1$.
 \el

 The following result was obtained in loc. cit (Lemma 3.5 and Corollary 3.9):

 \bp{prop-union}
 Let $W$ be any weak $V$-module. Set
 \begin{eqnarray}
 \Omega_{\infty}(W)=\bigcup_{n\in \N}\Omega_n(W)\subset W.
 \end{eqnarray}
 Then $\Omega_{\infty}(W)$ is a $V$-submodule of $W$. Furthermore, for any homogeneous vector $v\in V$ 
 and for any integer $k\ge \wt v$, $v_k$ is locally nilpotent on $\Omega_{\infty}(W)$.
 In particular,  $L(1)$ is locally nilpotent on $\Omega_{\infty}(W)$.
 \ep
 
 The following is straightforward: 
 
 \bp{filter-omega-n}
 Let $W$ be a weak $V$-module such that $W=\Omega_{\infty}(W)$. Form an $\N$-graded vector space
 \begin{eqnarray}
 {\rm gr}(W)=\bigoplus_{n\in \N}\left(\Omega_{n}(W)/\Omega_{n-1}(W)\right).
 \end{eqnarray}
 Then ${\rm gr}(W)$ is naturally an $\N$-graded weak $V$-module, and for $n\in \Z$, 
 $$\Omega_n({\rm gr}(W))=\bigoplus_{k\le n}\left(\Omega_{k}(W)/\Omega_{k-1}(W)\right).$$
 \ep
 
\br{rem-O-nW-trivial}
{\em  Note that if $W$ is an irreducible weak $V$-module such that $\Omega_n(W)\ne 0$ for some $n\in \N$, then $W=\Omega_{\infty}(W)$.
 Equivalently,  if $W$ is an irreducible weak $V$-module such that $W\ne \Omega_{\infty}(W)$, then $\Omega_n(W)= 0$ for all $n\in \N$.}
 \er

The following definition is due to Dong and Jiang (see \cite{dj1}):

\bd{def-O'-m-n}
{\em Let $W$ be a weak $V$-module and let $m,n\in \N$. 
Denote by $O^{\dagger}_{n,m}(W)$ the subspace of $W$, linearly spanned by vectors
\begin{eqnarray}
v\circ_{m}^{n}w:=\Res_{z}\frac{(1+z)^{m}}{z^{m+n+2}}Y((1+z)^{L(0)}v,z)w
\end{eqnarray}
for $v\in V,\ w\in W$. Furthermore, set
\begin{eqnarray}
O'_{n,m}(W)=O^{\dagger}_{n,m}(W)+(L(-1)+L(0)+m-n)W.
\end{eqnarray}}
\ed

Using Zhu's argument in \cite{zhu1} (cf. \cite{dlm-anv}, \cite{dj1}) we get
\begin{eqnarray}\label{general-O-n-m}
\Res_{z}\frac{(1+z)^{\wt v+m+s}}{z^{m+n+2+k}}Y(v,z)w\in O^{\dagger}_{n,m}(W)
\end{eqnarray}
for $v\in V,\ w\in W$ with $v$ homogeneous and for $s,k\in \N$ with $s\le k$.

\br{rem-inclusion}
{\em Let $m,n,p,q\in \N$. Then 
\begin{eqnarray}
O^{\dagger}_{q,m}(W)\subset O^{\dagger}_{n,m}(W) && \text{ if }q\ge n,\\
O^{\dagger}_{n,p}(W)\subset O^{\dagger}_{n,m}(W) && \text{ if } p\ge m.
\end{eqnarray}
On the other hand, we have 
\begin{eqnarray}
O_{n}(W)=O'_{n,n}(W).
\end{eqnarray}}
\er

Furthermore, we form two vector spaces, which will play a key role in this paper.

\bd{def-A-dagger-prime}
Let $W$ be a weak $V$-module and let $m,n\in \N$.  Define vector spaces
\begin{eqnarray}
&&A^{\dagger}_{n,m}(W)=W/O^{\dagger}_{n,m}(W),\\
&&A^{\diamond}_{n,m}(W)=W/O'_{n,m}(W).
\end{eqnarray}
\ed

\br{rem-DJ}  
{\em  In \cite{dj1}, Dong and Jiang introduced two more subspaces $O''_{n,m}(V)$, $O'''_{n,m}(V)$ of $V$, and defined
\begin{eqnarray*}
&&O_{n,m}(V)=O'_{n,m}(V)+O''_{n,m}(V)+O'''_{n,m}(V),\\
&&A_{n,m}(V):=V/O_{n,m}(V).
\end{eqnarray*}
It was proved therein that $A_{n,m}(V)$ has an $A_n(V)$-$A_m(V)$ bimodule structure in the sense that $A_{n,m}(V)$
has a left $A_{n}(V)$-module structure together with a commuting right $A_{m}(V)$-module structure. 
It was conjectured that $O_{n,m}(V)=O'_{n,m}(V)$, where the conjecture was confirmed in case $m=n$.
Since we use a different right $A_m(V)$-module action in this paper,  we skip the detailed definitions 
and precise theorems. }
\er

To end this section,  we formulate a linear algebra duality result which we shall use. 
Let $(A,*)$ be a non-associative algebra with an order-$2$ bijective linear operator $\theta$.
 Assume that $W$ is a vector space equipped with a bilinear map $A\times W\rightarrow W;\ (a,u)\mapsto au$.
Define a bilinear map
$$A\times W^*\rightarrow W^*; \quad (a,f)\mapsto af$$  by 
$$\< af,w\>=\<f,\theta(a)w\>\   \   \   \   \mbox{ for }a\in A, \ f\in W^{*},\  w\in W.$$

\bl{l-classical}
Let $W_0$ be a subspace of $W$ and view $(W/W_0)^{*}$ as a subspace of $W^{*}$ naturally through the quotient map
$W\rightarrow W/W_0$. Then $(W/W_0)^{*}$ is an $A$-stable subspace of $W^{*}$ if and only if $W_0$ is an $A$-stable subspace of $W$.
  \el
 
 Furthermore, assume that $J$ is an ideal of $A$ with $\theta(J)=J$ such that $A/J$ is an associative algebra and 
 $\theta$ is an algebra anti-automorphism. We have:
 
 \bl{l-classical-2}
Let $W_0$ be a subspace of $W$ such that $(W/W_0)^{*}$ is an $A$-stable subspace of $W^{*}$ with $J\cdot (W/W_0)^{*}=0$.
Then $J\cdot (W/W_0)=0$. Furthermore, $(W/W_0)^{*}$ is a module for $A/J$ viewed as an associative algebra 
if and only if $W/W_0$ is an $A/J$-module.
\el

\bl{l-classical-3}
Let $U$ be a vector space, let $\psi\in \End U$, and let  $U_1$, $U_2$ be subspaces of $U$. 
Denote by $\psi^{*}$ the dual operator of $\psi$ on $U^{*}$.
Then $\psi^{*}((U/U_1)^{*})\subset (U/U_2)^{*}$ if and only if $\psi(U_2)\subset U_1$.
\el

\section{The regular representation on $\mathfrak{D}_{(-1)}(W)$ and $A^{\dagger}_{n,m}(W)$}

In this section, we first recall the basic results on regular representations from \cite{Li-reg}, and we then 
introduce the level-$(m,n)$ vacuum space $\Omega_{m,n}(\mathfrak{D}_{(-1)}(W))$ of the regular representation 
module $\mathfrak{D}_{(-1)}(W)$ and give a canonical identification of $(A^{\dagger}_{n,m}(W))^{*}$ with 
$\Omega_{m,n}(\mathfrak{D}_{(-1)}(W))$. By making use this identification, we recover some of Dong-Jiang's results in \cite{dj1}.

Let $V$ be a vertex operator algebra and let $(W,Y_{W})$ be a weak $V$-module, 
both of which are fixed throughout  this section.
Define a linear map 
$$Y_{W}^{*}(\cdot,x):\  V\rightarrow (\End W^{*})[[x,x^{-1}]]$$ 
as in \cite{fhl} by the condition
\begin{eqnarray}
\< Y_{W}^{*}(v,x)f,w\>=\<f,Y_{W}(e^{xL(1)}(-x^{-2})^{L(0)}v,x^{-1})w\>
\end{eqnarray}
for $v\in V,\  f\in W^{*}, \ w\in W$.

\bd{dDz0}
{\em Let $z\in \C^{\times}$. 
Define $\mathfrak{D}_{(z)}(W)$ to consist of every $f\in W^{*}$, satisfying the condition that  for any $v\in V$,
there exists a nonnegative integer $k$ such that
\begin{eqnarray}\label{erational-prop}
(x-z)^{k}Y_{W}^{*}(v,x)f\in W^{*}((x)),
\end{eqnarray}
i.e.,  for any $v\in V$,
there exist nonnegative integers $k$ and $l$ such that
\begin{eqnarray}\label{erational-prop-2}
x^l(x-z)^{k}Y_{W}^{*}(v,x)f\in W^{*}[[x]].
\end{eqnarray}}
\ed

Let $v\in V,\ f\in \mathfrak{D}_{(z)}(W)$. Define
\begin{eqnarray}\label{eYR}
Y_{(z)}^{R}(v,x)f=(-z+x)^{-k}\left[(x-z)^{k}Y_{W}^{*}(v,x)f\right]\in W^{*}((x)),
\end{eqnarray}
where $k$ is any nonnegative integer such that (\ref{erational-prop}) holds and where as a convention
$$(-z+x)^{-k}=\sum_{j\ge 0}\binom{-k}{j}(-z)^{-k-j}x^{j}\in \C[[x]].$$
Note that (\ref{erational-prop}) implies that there is a nonnegative integer $l$ such that
\begin{eqnarray}\label{erational-prop-z0}
(x+z)^{l}Y_{W}^{*}(v,x+z)f\in W^{*}((x)).
\end{eqnarray}
(In fact, it was proved that (\ref{erational-prop}) and (\ref{erational-prop-z0}) are equivalent.)
Then define 
\begin{eqnarray}\label{eYL-z0}
Y_{(z)}^{L}(v,x)f=(z+x)^{-l}\left[(x+z)^{l}Y_{W}^{*}(v,x+z)f\right]\in W^{*}((x)).
\end{eqnarray}
The following is a key result in \cite{Li-reg}:

\bt{thm-reg}
Let $z\in \C^{\times}$.  Then both pairs $(\mathfrak{D}_{(z)}(W),Y_{(z)}^{L})$ and $(\mathfrak{D}_{(z)}(W),Y_{(z)}^{R})$ 
carry the structures of a weak $V$-module. Furthermore, the actions of $V$ under $Y_{(z)}^{R}(\cdot,x)$ and $Y_{(z)}^{L}(\cdot,x)$ 
commute, and $\mathfrak{D}_{(z)}(W)$ is naturally a weak module for the tensor product vertex operator algebra $V\otimes V$. 
\et

The following was also obtained therein (Lemma 3.23): 

\bl{tri-relation}
Let $v\in V,\ \alpha\in \mathfrak{D}_{(z)}(W)$. Then 
\begin{eqnarray}\label{tri-relation-alpha}
&&x_0^{-1}\delta\left(\frac{x-z}{x_0}\right)Y_W^{*}(v,x)\alpha-x_0^{-1}\delta\left(\frac{z-x}{-x_0}\right)Y^{R}_{(z)}(v,x)\alpha
\nonumber\\
&&\hspace{2cm} =z^{-1}\delta\left(\frac{x-x_0}{z}\right)Y^{L}_{(z)}(v,x_0)\alpha.
\end{eqnarray}
\el

For simplicity, we shall simply write $Y^R$ and $Y^L$ for $Y^R_{(z)}$ and $Y^L_{(z)}$ 
whenever it is clear from the context.

\bd{def-omega-m-n-M}
{\em Let $M$ be a weak $V\otimes V$-module. For $v\in V$, set
$$Y^{1}(v,x)=Y(v\otimes {\bf 1},x),\    \    \    \  Y^{2}(v,x)=Y({\bf 1}\otimes v,x).$$
For $m,n\in \N$, let $\Omega_{m,n}(M)$ consist of all vectors $w\in M$ such that
\begin{eqnarray}
x^{\wt v+m}Y^{1}(v,x)w,\quad   x^{\wt v+n}Y^{2}(v,x)w \in M[[x]]
\end{eqnarray}
for every homogeneous vector $v\in V$.}
\ed

From definition, we have 
\begin{eqnarray}
\Omega_{m,n}(M)=\Omega_m(M, Y^{1}(\cdot,x))\cap \Omega_n(M, Y^{2}(\cdot,x)).
\end{eqnarray}
As the actions of $V$ under $Y^{1}(\cdot,x)$ and $Y^{2}(\cdot,x)$ on $M$ commute, 
from Proposition \ref{prop-Omega-n-W}, $\Omega_{m,n}(M)$ is naturally an $A_{m}(V)$-module and an $A_{n}(V)$-module.
Furthermore,  $\Omega_{m,n}(M)$ is a (left) $A_{m}(V)\otimes A_{n}(V)$-module.

As the first result of this paper, we have:

\bp{p-A-D-space}
Let $W$ be a weak $V$-module and let $m,n\in \N$.  Set 
\begin{eqnarray}
A^{\dagger}_{n,m}(W)=W/O^{\dagger}_{n,m}(W),
\end{eqnarray}
a vector space. Let $f\in W^{*}$. Then $f\in (A^{\dagger}_{n,m}(W))^{*}$ if and only if
\begin{eqnarray}\label{2.16}
x^{\wt v+n}(x+1)^{\wt v+m}Y_{W}^{*}(v,x)f\in W^{*}[[x]]
\end{eqnarray}
for every homogeneous vector $v\in V$. Furthermore, we have
\begin{eqnarray}
(A^{\dagger}_{n,m}(W))^{*}= \Omega_{m,n}(\mathfrak{D}_{(-1)}(W))
\end{eqnarray}
as subspaces of $W^{*}$.
\ep

\begin{proof} Suppose $f\in (A^{\dagger}_{n,m}(W))^{*}\subset W^{*}$, i.e.,  $f\in W^{*}$ such that $f|_{O^{\dagger}_{n,m}(W)}=0$.
Let $v\in V$ be homogeneous. For every $w\in W$ and for any $k\in \N$, using (\ref{general-O-n-m}) we get
\begin{eqnarray*}
&&\Res_{x}x^{\wt v+n+k}(x+1)^{\wt v+m}\<Y_{W}^{*}(v,x)f,w\>\nonumber\\
&=&\Res_{x}x^{\wt v+n+k}(x+1)^{\wt v+m}\<f,Y_{W}(e^{xL(1)}(-x^{-2})^{L(0)}v,x^{-1})w\>\nonumber\\
&=&(-1)^{\wt v}\Res_{x}x^{n+k-\wt v} (x+1)^{\wt v+m} 
\<f,Y_{W}(e^{xL(1)}v,x^{-1})w\>\nonumber\\
&=&(-1)^{\wt v}\Res_{z}\frac{(1+z)^{\wt v+m}}{z^{m+n+k+2}} 
\<f,Y_{W}(e^{z^{-1}L(1)}v,z)w\>\nonumber\\
&=&(-1)^{\wt v}\sum_{i\ge 0}\frac{1}{i!}\Res_{z}\frac{(1+z)^{\wt L(1)^{i}v+m+i}}{z^{m+n+k+2+i}} 
\<f,Y_{W}(L(1)^{i}v,z)w\>\nonumber\\
&=&0,
\end{eqnarray*}
which implies (\ref{2.16}).
Conversely, assume that $f\in W^{*}$ such that (\ref{2.16}) holds for every homogeneous vector $v\in V$. 
Let $u\in V$ be homogeneous. Then for any $w\in W$ we have
\begin{eqnarray*}
&&\Res_z\frac{(1+z)^{\wt u +m}}{z^{m+n+2}} \<f, Y_W(u,z)w\>\nonumber\\
&=&\Res_z\frac{(1+z)^{\wt u +m}}{z^{m+n+2}}\<Y_W^{*}(e^{zL(1)}(-z^{-2})^{L(0)}u,z^{-1})f,w\>\nonumber\\
&=&(-1)^{\wt u}\Res_x x^{\wt u +n}(x+1)^{\wt u +m} \<Y_W^{*}(e^{x^{-1}L(1)}u,x)f,w\>\nonumber\\
&=&(-1)^{\wt u}\sum_{i\ge 0}\frac{1}{i!}\Res_xx^{\wt u+n-i}(x+1)^{\wt u+m}\<Y_W^{*}(L(1)^{i}u,x)f,w\>\nonumber\\
&=&(-1)^{\wt u}\sum_{i\ge 0}\frac{1}{i!}\Res_xx^{\wt L(1)^iu+n}(x+1)^{\wt L(1)^iu+m+i}\<Y_W^{*}(L(1)^{i}u,x)f,w\>\nonumber\\
&=&0.
\end{eqnarray*}
This proves $f|_{O^{\dagger}_{n,m}(W)}=0$. That is, $f\in (A^{\dagger}_{n,m}(W))^{*}$. This proves the first assertion.

Now, let $f\in W^{*}$ such that (\ref{2.16}) holds for every homogeneous vector $v\in V$.
Note that (\ref{2.16}) implies $f\in \mathfrak{D}_{(-1)}(W)$. Furthermore, from the definition of $Y^{R}(v,x)f$ we have
\begin{eqnarray}\label{eYR=Y*}
x^{\wt v+n}(1+x)^{\wt v+m}Y^{R}(v,x)f=x^{\wt v+n}(x+1)^{\wt v+m}Y_{W}^{*}(v,x)f.
\end{eqnarray}
Combining this with  (\ref{2.16}) we get
$$x^{\wt v+n}(1+x)^{\wt v+m}Y^{R}(v,x)f\in \mathfrak{D}_{(-1)}(W)[[x]],$$
which implies 
\begin{eqnarray}\label{YRcondition}
x^{\wt v+n}Y^{R}(v,x)f\in \mathfrak{D}_{(-1)}(W)[[x]].
\end{eqnarray}
On the other hand, from (\ref{2.16}) we get
\begin{eqnarray*}
(x-1)^{\wt v+n}x^{\wt v+m}Y_{W}^{*}(v,x-1)f\in W^{*}[[x]].
\end{eqnarray*}
From the definition of $Y^{L}(v,x)f$ we have
\begin{eqnarray}\label{YL=Y*}
(-1+x)^{\wt v+n}x^{\wt v+m}Y^{L}(v,x)f=(x-1)^{\wt v+n}x^{\wt v+m}Y_{W}^{*}(v,x-1)f.  
\end{eqnarray}
By the same reasoning, we get
\begin{eqnarray}\label{YLcondition}
x^{\wt v+m}Y^{L}(v,x)f\in \mathfrak{D}_{(-1)}(W)[[x]].
\end{eqnarray}
With (\ref{YRcondition}) and  (\ref{YLcondition}) we conclude $f\in \Omega_{m,n}(\mathfrak{D}_{(-1)}(W))$.
This proves $(A^{\dagger}_{n,m}(W))^{*}\subset \Omega_{m,n}(\mathfrak{D}_{(-1)}(W))$.
On the other hand, let $f\in  \Omega_{m,n}(\mathfrak{D}_{(-1)}(W))$. For any homogeneous vector $v\in V$, 
using  (\ref{tri-relation-alpha}), (\ref{YRcondition}) and  (\ref{YLcondition}) we see that (\ref{2.16}) holds
(multiplying both sides by $x^{\wt v+n}x_0^{\wt v+m}$ and then applying $\Res_{x_0}$). Then $f\in (A^{\dagger}_{n,m}(W))^{*}$.
This proves $\Omega_{m,n}(\mathfrak{D}_{(-1)}(W))\subset (A^{\dagger}_{n,m}(W))^{*}$, concluding the proof.
\end{proof}

 Next, we refine Proposition \ref{p-A-D-space}, to get an identification of $A_{m}(V)\otimes A_n(V)$-modules. 
 To this end, we shall need the following two results from \cite{Li-Anv} (cf. \cite{Li-reg},  Remark 2.10): 
 
 \bl{lemma-e-L(1)}
 Let $(E,Y_E)$ be a weak $V$-module and let $z_0\in \C$. For $v\in V$, define
 \begin{eqnarray}\label{deformed-module}
 Y_E^{[z_0]}(v,x)=Y_E(e^{-z_0(1+z_0x)L(1)}(1+z_0x)^{-2L(0)}v,x/(1+z_0x)).
 \end{eqnarray}
 Then the pair $(E, Y_E^{[z_0]})$ carries the structure of a weak $V$-module,  and  for homogeneous vector $v\in V$ and for $m\in \Z$,
 \begin{eqnarray}\label{degree-0-formula}
\Res_x x^mY_E^{[z_0]}(v,x)=\Res_x x^m(1-z_0x)^{2\wt v-m-2}Y_E(e^{-z_0(1-z_0x)^{-1}L(1)}v,x). 
 \end{eqnarray}
 Furthermore, if $L(1)$ is locally nilpotent on $E$, $e^{-z_0L(1)}$ is a $V$-module isomorphism from $(E,Y_E)$ to $(E,Y_E^{[z_0]})$.
 \el
 
 \bp{p-anv-omega}
 Let $(E,Y_E)$ be a weak $V$-module and let $z_0\in \C$. Then for every $n\in \N$,
 \begin{eqnarray}
 \Omega_n(E,Y_E)= \Omega_n(E,Y_E^{[z_0]}).
 \end{eqnarray}
 Furthermore, $e^{-z_0L(1)}$ is an $A_n(V)$-module isomorphism from $\Omega_n(E,Y_E)$ to $\Omega_n(E,Y_E^{[z_0]})$.
 \ep
 
 
 As an immediate consequence of Lemma \ref{lemma-e-L(1)} and Proposition \ref{p-anv-omega}, we have:
 
\bc{z-0-anv-module}
Let $(E,Y_E)$ be a weak $V$-module and let $z_0\in \C$. 
Define a bilinear map from $V\times E$ to $E$ by
\begin{eqnarray}
v\bullet_{(z_0)} w=\Res_x x^{\wt v-1}(1-z_0x)^{\wt v-1}Y_E(e^{-z_0(1-z_0x)^{-1}L(1)}v,x)w 
\end{eqnarray}
for $v\in V,\ w\in E$ with $v$ homogeneous. Then this bilinear map gives rise to an $A_n(V)$-module
structure on $\Omega_n(E)$ for every $n\in \N$. We denote this $A_n(V)$-module by $\Omega_n^{[z_0]}(E)$.
 \ec
  
 \bd{Am-action}
{\em Define a bilinear map $\bar{*}_{m,n}: V\times W\rightarrow W; \ (v,w)\mapsto v\bar{*}_{m,n}w$ by
\begin{eqnarray} 
v\bar{*}_{m,n}w=\Res_x \sum_{i=0}^{m}\binom{-n-1}{i}
(-1)^{n+i}\frac{(1+x)^{i-1}}{x^{n+i+1}}Y_W\left((1+x)^{L(0)}v,x\right)w
\end{eqnarray}
for $v\in V,\ w\in W$.}
\ed

\bl{Y-L-DJ}
Let $v\in V,\ f\in A^{\dagger}_{n,m}(W)^{*}\subset W^*,\ w\in W$ with $v$ homogeneous. Then
\begin{eqnarray}
\Res_{x}x^{\wt v-1}\<(Y^{L})^{[1]}(v,x)f,w\>=\<f, v\bar{*}_{m,n}w\>.
\end{eqnarray}
\el
 
\begin{proof} Note that $x^{\wt u+m}Y^{L}(u,x)f\in \mathfrak{D}_{(-1)}(W)[[x]]$ for any homogeneous vector $u\in V$.
As $(-1+x)^{p}\in \C[[x]]$ and $\wt (L(1)^rv)=\wt v-r$ for any $p\in \Z,\ r\in \N$, we have
$$\Res_xx^{\wt v+j-1} (-1+x)^{\wt v+n}Y^{L}(e^{(-1+x)^{-1}L(1)}v,x)f=0 $$ 
 for all $j>m$. From (\ref{YL=Y*}), we have
 $$(-1+x)^{\wt u+n}Y^{L}(u,x)f=(x-1)^{\wt u+n}Y_W^{*}(u,x-1)f$$
 for any homogeneous vector $u\in V$, which yields
 $$(-1+x)^{\wt v+n}Y^{L}(e^{(-1+x)^{-1}L(1)}v,x)f=(x-1)^{\wt v+n}Y_W^{*}(e^{(x-1)^{-1}L(1)}v,x-1)f.$$
 Then using (\ref{degree-0-formula}) we get
 \begin{eqnarray*}
&&\Res_{x}x^{\wt v-1}\<(Y^{L})^{[1]}(v,x)f,w\>\\
&=&\Res_x (-1)^{\wt v-1}x^{\wt v-1} (-1+x)^{\wt v-1}\<Y^{L}(e^{(-1+x)^{-1}L(1)}v,x)f,w\>\\
&=&\Res_x \sum_{i=0}^{\infty}\binom{-n-1}{i}
(-1)^{\wt v+n-i}x^{\wt v+i-1} (-1+x)^{\wt v+n}\<Y^{L}(e^{(-1+x)^{-1}L(1)}v,x)f,w\>\\
&=&\Res_x \sum_{i=0}^{m}\binom{-n-1}{i}
(-1)^{\wt v+n-i}x^{\wt v+i-1} (-1+x)^{\wt v+n}\<Y^{L}(e^{(-1+x)^{-1}L(1)}v,x)f,w\>\\
&=&\Res_x \sum_{i=0}^{m}\binom{-n-1}{i}
(-1)^{\wt v+n-i}x^{\wt v+i-1} (x-1)^{\wt v+n}\<Y_W^{*}(e^{(x-1)^{-1}L(1)}v,x-1)f,w\>\\
&=&\Res_x \sum_{i=0}^{m}\binom{-n-1}{i}
(-1)^{\wt v+n-i}(x+1)^{\wt v+i-1} x^{\wt v+n}\<Y_W^{*}(e^{x^{-1}L(1)}v,x)f,w\>\\
&=&\Res_x \sum_{i=0}^{m}\binom{-n-1}{i}
(-1)^{\wt v+n-i}(x+1)^{\wt v+i-1} x^{\wt v+n}\\
&&\quad\quad \times \<f, Y_W(e^{xL(1)}(-x^{-2})^{L(0)}e^{x^{-1}L(1)}v,x^{-1})w\>\\
&=&\Res_x \sum_{i=0}^{m}\binom{-n-1}{i}
(-1)^{\wt v+n-i}(x+1)^{\wt v+i-1} x^{\wt v+n}\<f, Y_W((-x^{2})^{-L(0)}v,x^{-1})w\>\\
&=&\Res_x \sum_{i=0}^{m}\binom{-n-1}{i}
(-1)^{n-i}(x+1)^{\wt v+i-1} x^{-\wt v+n}\<f, Y_W(v,x^{-1})w\>\\
&=&\Res_x \sum_{i=0}^{m}\binom{-n-1}{i}
(-1)^{n-i}(x^{-1}+1)^{\wt v+i-1} x^{\wt v-n-2}\<f, Y_W(v,x)w\>\\
&=&\Res_x \sum_{i=0}^{m}\binom{-n-1}{i}
(-1)^{n+i}\frac{(1+x)^{\wt v+i-1}}{x^{n+i+1}} \<f, Y_W(v,x)w\>\\
&=&\<f, v\bar{*}_{m,n}w\>,
 \end{eqnarray*}
 as desired.
 \end{proof}
 
 Combining Lemma \ref{Y-L-DJ} with Lemmas \ref{l-classical} and  \ref{l-classical-2}, we immediately obtain:
 
 \bp{cor-L-connection}
The bilinear map $V\times W\rightarrow W;\  (v,w)\mapsto v\cdot w:=v\bar{*}_{m,n}w$ 
reduces to a bilinear map from $A_m(V)\times A^{\dagger}_{n,m}(W)$ to $A^{\dagger}_{n,m}(W)$ such that 
$A^{\dagger}_{n,m}(W)$ is a right $A_m(V)$-module.
 \ep
 
 Recall the following bilinear map from \cite{dj1}:

\bd{An-action}
{\em Define a bilinear map $\bar{*}_{m}^{n}: V\times W\rightarrow W; \ (v,w)\mapsto v\bar{*}_{m}^{n}w$ by
\begin{eqnarray}
v\bar{*}_{m}^{n}w=\sum_{i=0}^{n}\binom{-m-1}{i}\Res_{z}\frac{(1+z)^{m}}{z^{m+i+1}}Y((1+z)^{L(0)}v,z)w
\end{eqnarray}
for $v\in V,\ w\in W$.}
\ed

The following is a connection between the product $v\bar{*}_{m,n}w$ 
and the product $v\bar{*}_{m}^nw$ for $v\in V,\ w\in W$:

\bl{left-right-diff}
For any homogeneous vector $v\in V$ and for any $w\in W$, we have
  \begin{eqnarray}
 v\bar{*}_{m,n}w=v\bar{*}_{m}^nw-\Res_x (1+x)^{\wt v-1}Y(v,x)w.
 \end{eqnarray}
\el

\begin{proof}  Recall the following identity from \cite{dj1} (Proposition 5.1):
 \begin{eqnarray}\label{DJ-formula}
\sum_{i=0}^m\binom{-n-1}{i}\frac{(1+x)^{n+1}}{x^{n+i+1}}-\sum_{i=0}^n\binom{-m-1}{i}(-1)^{m+i}\frac{(1+x)^i}{x^{m+i+1}}=1.
 \end{eqnarray}
Using this (with $m$ and $n$ switched) we have
\begin{eqnarray} 
&& \sum_{i=0}^m\binom{-n-1}{i}(-1)^{n+i}\frac{(1+x)^{\wt v+i-1}}{x^{n+i+1}}\nonumber\\
&=&(1+x)^{\wt v-1}\left[-1+\sum_{i=0}^n\binom{-m-1}{i}\frac{(1+x)^{m+1}}{x^{m+i+1}}\right]\nonumber\\
&=&-(1+x)^{\wt v-1}+\sum_{i=0}^n\binom{-m-1}{i}\frac{(1+x)^{\wt v+m}}{x^{m+i+1}}.
\end{eqnarray}
Then by the definitions of $v\bar{*}_{m,n}w$ and $v\bar{*}_{m}^nw$ we obtain
 \begin{eqnarray*}
 v\bar{*}_{m,n}w=v\bar{*}_{m}^nw-\Res_x (1+x)^{\wt v-1}Y(v,x)w.
 \end{eqnarray*}
 as desired.
 \end{proof}

We have the following analogue of Lemma \ref{Y-L-DJ}:

 \bl{R-dj-relation}
 Let $v\in V,\ f\in A^{\dagger}_{n,m}(W)^{*}\subset W^*,\ w\in W$ with $v$ homogeneous. Then
  \begin{eqnarray}\label{R-dj-relation-formula}
\Res_{x}x^{\wt v-1}\<(Y^{R})^{[-1]}(v,x)f,w\>= \<f, \theta(v)\bar{*}_{m}^{n}w\>.
\end{eqnarray}
\el

\begin{proof} Recall that for any homogeneous vector $u\in V$, we have 
\begin{eqnarray*}
&&\quad \quad  x^{\wt u+n}Y^R(u,x)f\in \mathfrak{D}_{(-1)}[[x]],\\ 
&&(1+x)^{\wt u+m}Y^{R}(u,x)f=(x+1)^{\wt u+m}Y^{*}(u,x)f.
\end{eqnarray*}
Note that $(1+x)^{p}\in \C[[x]]$ for $p\in \Z$ and $\wt (L(1)^rv)=\wt v-r$ for $r\in \N$. 
Using all of these facts and (\ref{degree-0-formula}), we obtain
  \begin{eqnarray*}
&&\Res_{x}x^{\wt v-1}\<(Y^{R})^{[-1]}(v,x)f,w\>\\
&=&\Res_x x^{\wt v-1} (1+x)^{\wt v-1}\<Y^{R}(e^{(1+x)^{-1}L(1)}v,x)f,w\>\\
&=&\Res_x \sum_{i=0}^{\infty}\binom{-m-1}{i}
x^{\wt v+i-1} (1+x)^{\wt v+m}\<Y^{R}(e^{(1+x)^{-1}L(1)}v,x)f,w\>\\
&=&\Res_x \sum_{i=0}^{n}\binom{-m-1}{i}
x^{\wt v+i-1} (1+x)^{\wt v+m}\<Y^{R}(e^{(1+x)^{-1}L(1)}v,x)f,w\>\\
&=&\Res_x \sum_{i=0}^{n}\binom{-m-1}{i}
x^{\wt v+i-1} (x+1)^{\wt v+m}\<Y_W^{*}(e^{(x+1)^{-1}L(1)}v,x)f,w\>\\
&=&\Res_x \sum_{i=0}^{n}\binom{-m-1}{i}
x^{\wt v+i-1} (x+1)^{\wt v+m}\\
&&\quad\quad \times \<f, Y_W(e^{xL(1)}(-x^{-2})^{L(0)}e^{(x+1)^{-1}L(1)}v,x^{-1})w\>\\
&=&\Res_x \sum_{i=0}^{n}\binom{-m-1}{i}
x^{\wt v+i-1} (x+1)^{\wt v+m}  \<f, Y_W(e^{x(x+1)^{-1}L(1)}  (-x^{-2})^{L(0)}v,x^{-1})w\>\\
&=&\Res_x \sum_{i=0}^{n}\binom{-m-1}{i}
 \frac{(1+x)^{\wt v+m}}{x^{m+i+1}} \<f, Y_W(e^{(1+x)^{-1}L(1)}(-1)^{L(0)}v,x)w\>\\
&=&\sum_{r\ge 0}\frac{1}{r!}\Res_x \sum_{i=0}^{n}\binom{-m-1}{i}
\frac{(1+x)^{\wt (L(1)^rv)+m}} {x^{m+i+1}}  \<f, Y_W(L(1)^r(-1)^{L(0)}v,x)w\>\\
&=& \<f, \theta(v)\bar{*}_{m}^{n}w\>,
\end{eqnarray*}
as desired.
 \end{proof}
 
By invoking the linear algebra duality again, we immediately recover the following result of \cite{dj1}:
 
\bp{left-module-bar}
The bilinear map $V\times W\rightarrow W; \  (v,w)\mapsto v\bar{*}_{m}^{n}w$ reduces to a bilinear map
from $A_n(V)\times A^{\dagger}_{n,m}(W)$ to $A^{\dagger}_{n,m}(W)$ such that
$A^{\dagger}_{n,m}(W)$ becomes a (left) $A_n(V)$-module.
\ep
 
As the actions of $V$ on $\mathfrak{D}_{(-1)}(W)$ under $Y^{L}$ and $Y^{R}$ commute, it follows from (\ref{deformed-module}) that
the actions of $V$ on $\mathfrak{D}_{(-1)}(W)$ under $(Y^{L})^{[1]}$ and $(Y^{R})^{[-1]}$ commute.
Then using Lemmas \ref{Y-L-DJ} and \ref{R-dj-relation} we immediately have:
 
\bp{prop-A-dagger}
Let $m,n\in \N$. Then on the space $A^{\dagger}_{n,m}(W)$ the left module action of $A_n(V)$ established in Proposition \ref{left-module-bar} 
and the right module action of $A_m(V)$ established in Proposition \ref{cor-L-connection} commute. 
 \ep

To summarize, as the main result of this section we have:

\bt{thm-Amn}
Let $W$ be a weak $V$-module and let $m,n\in \N$. 
Then $\Omega_{m,n}(\mathfrak{D}_{(-1)}(W))$ is an $A_{m}(V)\otimes A_{n}(V)$-module 
with the actions of $A_m(V)$ and $A_n(V)$ given respectively by 
\begin{eqnarray}
&&v\bullet_L f=\Res_x x^{\wt v-1}(1-x)^{\wt v-1}Y^{L}(e^{(-1+x)^{-1}L(1)}v,x)f,\label{vLf}\\
&&v\bullet_R f=\Res_x x^{\wt v-1}(1+x)^{\wt v-1}Y^{R}(e^{(1+x)^{-1}L(1)}v,x)f\label{vRf}
 \end{eqnarray}
for $v\in V,\  f\in \mathfrak{D}_{(-1)}(W)$. 
Denote this $A_{m}(V)\otimes A_{n}(V)$-module by $\Omega_{m,n}^{[1,-1]}(\mathfrak{D}_{(-1)}(W))$.
On the other hand, the subspace $A^{\dagger}_{n,m}(W)^{*}$ of $W^{*}$ 
is an $A_{m}(V)\otimes A_{n}(V)$-module with the actions of $A_m(V)$ and $A_n(V)$ given respectively by
\begin{eqnarray}
&&\<v\cdot_L f,w\>=\< f, v\bar{*}_{m,n}w\>,\label{v-dot-L}\\
&&\<v\cdot_R f,w\>=\< f, \theta(v)\bar{*}_{m}^{n}w\>\label{v-dot-R}
\end{eqnarray}
for $v\in V,\ f\in A^{\dagger}_{n,m}(W)^{*},\  w\in W$.
Furthermore, we have
\begin{eqnarray}
\Omega_{m,n}^{[1,-1]}(\mathfrak{D}_{(-1)}(W))=(A^{\dagger}_{n,m}(W))^{*}
\end{eqnarray} 
as $A_{m}(V)\otimes A_{n}(V)$-modules. 
\et

\begin{proof}  From definition, we have
\begin{eqnarray*}
\Omega_{m,n}(\mathfrak{D}_{(-1)}(W))
=\Omega_{m}(\mathfrak{D}_{(-1)}(W),Y^L)\cap \Omega_{n}(\mathfrak{D}_{(-1)}(W),Y^R)
\end{eqnarray*}
as vector spaces. In view of Corollary \ref{z-0-anv-module}, 
$\Omega_{m}(\mathfrak{D}_{(-1)}(W),Y^L)$ is an $A_m(V)$-module with the action given by
(\ref{vLf}), and $\Omega_{n}(\mathfrak{D}_{(-1)}(W),Y^R)$ is an $A_n(V)$-module with the action given by (\ref{vRf}).
Recall that the actions of $V$ under $Y^R$ and $Y^L$ on $\mathfrak{D}_{(-1)}(W)$ commute.
Then it follows that $\Omega_{m,n}(\mathfrak{D}_{(-1)}(W))$ is an $A_m(V)\otimes A_n(V)$-module 
with the actions of $A_m(V)$ and $A_n(V)$ given respectively by (\ref{vLf}) and (\ref{vRf}). This proves the first assertion.

As $\theta$ gives rise to an anti-automorphism of $A_n(V)$, with Proposition \ref{prop-A-dagger}, we conclude that
$A^{\dagger}_{n,m}(W)^{*}$ is an $A_m(V)\otimes A_n(V)$-module with the actions of $A_m(V)$ and $A_n(V)$ given respectively by
(\ref{v-dot-L}) and (\ref{v-dot-R}). The furthermore assertion follows from Propositions \ref{Y-L-DJ}  and \ref{R-dj-relation}.
\end{proof}

\br{Omega-direct-limit}
{\em Let $W$ be a weak $V$-module. Set 
\begin{eqnarray}
\Omega_{\infty}(\mathfrak{D}_{(-1)}(W))=\cup_{m,n\in \N}\Omega_{m,n}(\mathfrak{D}_{(-1)}(W)),
\end{eqnarray}
which is a $V\otimes V$-submodule of $\mathfrak{D}_{(-1)}(W)$.
Equip $\N\times \N$ with the partial order $\le$, where $(m,n)\le (p,q)$ if and only if $m\le p$ and $n\le q$.
 Then $\Omega_{m,n}(\mathfrak{D}_{(-1)}(W))$ for $m,n\in \N$ form a direct system over $\N\times \N$,
 for which $\Omega_{\infty}(\mathfrak{D}_{(-1)}(W))$ is a direct limit. }
 \er

\br{rem-A(W)}
{\em For $n\in \N$, as $O^{\dagger}_n(W)=O^{\dagger}_{n,n}(W)$, we have $A^{\dagger}_n(W)=A^{\dagger}_{n,n}(W)$. 
Then
\begin{eqnarray}
(A^{\dagger}_n(W))^{*}=(A^{\dagger}_{n,n}(W))^{*}=\Omega_{n,n}(\mathfrak{D}_{(-1)}(W)).
\end{eqnarray}
Consequently, we have
\begin{eqnarray}
\Omega_{\infty}(\mathfrak{D}_{(-1)}(W))=\cup_{n\in \N}(A^{\dagger}_n(W))^{*}.
\end{eqnarray}}
\er

Recall the weak $V$-module structure $(Y^R)^{[-1]}(\cdot,x)$ on $\mathfrak{D}_{(-1)}(W)$, where 
\begin{eqnarray}
\Res_xx^k(Y^R)^{[-1]}(v,x)=\Res_xx^k(1+x)^{2\wt v-k-2}Y^R(e^{(1+x)^{-1}L(1)}v,x)
\end{eqnarray}
for homogeneous vector $v\in V$ and for $k\in \Z$. For any $u\in V$, write
\begin{eqnarray}
(Y^R)^{[-1]}(u,x)=\sum_{k\in\Z}u^{R[-1]}_kx^{-k-1}.
\end{eqnarray}

Using Lemma \ref{lemma-omega-union} and Proposition \ref{p-anv-omega} we immediately have:

\bl{lem-YR[-1]}
Let $m,n\in \N,\ p\in \Z$. Then for any homogeneous vector $v\in V$,
\begin{eqnarray}
v_{\wt v-1+p}^{R[-1]}\cdot \Omega_{m,n}(\mathfrak{D}_{(-1)}(W))\subset \Omega_{m,n-p}(\mathfrak{D}_{(-1)}(W)),
\end{eqnarray}
where  by convention $\Omega_{m,k}(\mathfrak{D}_{(-1)}(W))=0$ for $k<0$.
\el

The following is a slight variation of Dong and Jiang's notion $u*_{m,p}^{n}v$ (see  \cite{dj1}):

\bd{def-general-product}
{\em Let $m,n\in \N,\ p\in \Z$. For $u\in V,\ w\in W$, define
\begin{eqnarray}\label{deform-dj}
u[p]\bar{*}_{m}^{n}w=\Res_x \sum_{i=0}^{n+|p|}\binom{-m-p-1}{i}
\frac{(1+x)^{m}} {x^{p+m+i+1}}Y_W((1+x)^{L(0)}u,x)w\in W.
\end{eqnarray}}
\ed

We have the following generalization of Lemma \ref{R-dj-relation}:

\bl{lem-general-YR}
Let $f\in \Omega_{m,n}(\mathfrak{D}_{(-1)}(W))=(A_{n,m}^{\dagger}(W))^{*}$ with $m,n\in \N$. 
Then for any $v\in V, \  p\in \Z,\ w\in W$, we have
\begin{eqnarray}
\Res_x x^{p-1}\<(Y^{R})^{[-1]}(x^{L(0)}v,x)f,w\>= \<f, \theta(v)[p]\bar{*}_{m}^{n}w\>.
\end{eqnarray}
\el

\begin{proof}  Let $v\in V,\ p\in \Z,\ w\in W$ with $v$ homogeneous. Note that
\begin{eqnarray*}
&&x^{\wt v-1+p+i} (1+x)^{\wt v+m}Y^{R}(e^{(1+x)^{-1}L(1)}v,x)f\\
&=&\sum_{r\ge 0}\frac{1}{r!}x^{\wt (L(1)^rv)-1+r+p+i} (1+x)^{\wt (L(1)^rv)+m}Y^{R}(L(1)^{r}v,x)f\\
&\in & W^{*}[[x]] 
\end{eqnarray*}
for $i>n+|p|$.
Then using (\ref{degree-0-formula}) we get
\begin{eqnarray*}
&&\Res_{x}x^{p-1}\<(Y^{R})^{[-1]}(x^{L(0)}v,x)f,w\>\\
&=&\Res_x x^{\wt v-1+p} (1+x)^{\wt v-1-p}\<Y^{R}(e^{(1+x)^{-1}L(1)}v,x)f,w\>\\
&=&\Res_x \sum_{i=0}^{\infty}\binom{-m-p-1}{i}
x^{\wt v-1+p+i} (1+x)^{\wt v+m}\<Y^{R}(e^{(1+x)^{-1}L(1)}v,x)f,w\>\\
&=&\Res_x \sum_{i=0}^{n+|p|}\binom{-m-p-1}{i}
x^{\wt v-1+p+i} (1+x)^{\wt v+m}\<Y^{R}(e^{(1+x)^{-1}L(1)}v,x)f,w\>\\
&=&\Res_x \sum_{i=0}^{n+|p|}\binom{-m-p-1}{i}
x^{p+\wt v-1+i} (x+1)^{\wt v+m}\<Y_W^{*}(e^{(x+1)^{-1}L(1)}v,x)f,w\>\\
&=&\Res_x \sum_{i=0}^{n+|p|}\binom{-m-p-1}{i}
x^{\wt v-1+p+i} (x+1)^{\wt v+m}\\
&&\quad\quad \times \<f, Y_W(e^{xL(1)}(-x^{-2})^{L(0)}e^{(x+1)^{-1}L(1)}v,x^{-1})w\>\\
&=&\Res_x \sum_{i=0}^{n+|p|}\binom{-m-p-1}{i}
x^{\wt v-1+p+i} (x+1)^{\wt v+m}  \<f, Y_W(e^{x(x+1)^{-1}L(1)}  (-x^{-2})^{L(0)}v,x^{-1})w\>\\
&=&\Res_x \sum_{i=0}^{n+|p|}\binom{-m-p-1}{i}
 \frac{(1+x)^{\wt v+m}}{x^{p+m+i+1}} \<f, Y_W(e^{(1+x)^{-1}L(1)}(-1)^{L(0)}v,x)w\>\\
&=&\sum_{r\ge 0}\frac{1}{r!}\Res_x \sum_{i=0}^{n+|p|}\binom{-m-p-1}{i}
\frac{(1+x)^{\wt (L(1)^rv)+m}} {x^{p+m+i+1}}  \<f, Y_W(L(1)^r(-1)^{L(0)}v,x)w\>\\
&=& \<f, \theta(v)[p]\bar{*}_{m}^{n}w\>,
\end{eqnarray*}
as desired.
\end{proof}

Note that for homogeneous vector $v\in V$, as $\theta^2=1$ we also have
\begin{eqnarray}\label{other-direction}
\<f, v[p]\bar{*}_{m}^{n}w\>&=&\Res_x x^{p-1}\< (Y^R)^{[-1]}(x^{L(0)}\theta(v),x)f,w\>\nonumber\\
&=&\sum_{r\ge 0}\frac{1}{r!}\< (L(1)^r(-1)^{L(0)}v)^{R[-1]}_{\wt v-r-1+p}f,w\>.
\end{eqnarray}
Using Lemmas \ref{lem-YR[-1]}, \ref{lem-general-YR} and \ref{l-classical-3} we immediately recover
the following result of \cite{dj1}: 

\bc{coro-O-m-n}
Let $m,n\in \N,\ p\in \Z$, and let $u\in V$. Then 
\begin{eqnarray}
u[p]\bar{*}_{m}^{n} O^{\dagger}_{n,m}(W)\subset O^{\dagger}_{n+p,m}(W),
\end{eqnarray}
where by definition $O^{\dagger}_{k,m}(W)=W$ for $k<0$  (so that $A_{k,m}^{\dagger}(W)^{*}=0=\Omega_{m,k}^{\dagger}(\mathfrak{D}_{(-1)}(W))$).
\ec

\section{$A^{\diamond}_{n,m}(W)$ and $\Omega^{\diamond}_{m,n}(\mathfrak{D}_{(-1)}(W))$}

Here, we continue to identify $(A^{\diamond}_{n,m}(W))^{*}$ with a particular subspace 
$\Omega^{\diamond}_{m,n}(\mathfrak{D}_{(-1)}(W))$ of $\Omega_{m,n}(\mathfrak{D}_{(-1)}(W))$.
For each fixed $m\in \N$, the sum of subspaces $\Omega^{\diamond}_{m,n}(\mathfrak{D}_{(-1)}(W))$ for $n\in \N$ is shown to be a direct sum and
an $\N$-graded weak $V$-module. Using this connection,  we obtain an
$\N$-graded weak $V$-module structure on $A^{\diamond}_{\Box,m}(W):=\bigoplus_{n\in \N}A^{\diamond}_{n,m}(W)$, which recovers a result of
\cite{dj1}. 

Let $(W,Y_W)$ be a weak $V$-module as before. Write
$$Y^{*}_W(\omega,x)=\sum_{k\in \Z}L^{*}(k)x^{-k-2}$$
on $W^{*}$, where $\omega$ is the conformal vector of $V$. We have
$$\< L^{*}(k)\alpha,u\>=\<\alpha, L(-k)u\>\quad \text{ for }\alpha\in W^{*},\ u\in W.$$
Recall the weak $V$-module structures  $(Y^L)^{[1]}(\cdot,x)$ and $(Y^R)^{[-1]}(\cdot,x)$ on $\mathfrak{D}_{(-1)}(W)$. Write
\begin{eqnarray}
(Y^L)^{[1]}(\omega,x)=\sum_{k\in \Z}L_l^{[1]}(k)x^{-k-2},\quad (Y^R)^{[-1]}(\omega,x)=\sum_{k\in \Z}L_r^{[-1]}(k)x^{-k-2}.
\end{eqnarray}
Recall that $O'_{n,m}(W)=O^{\dagger}_{n,m}(W)+(L(-1)+L(0)+m-n)W$ for $m,n\in \N$. We have:

\bp{p-L(-1)+L(0)}
The following relation holds for any $f\in \mathfrak{D}_{(-1)}(W)$:
\begin{eqnarray}\label{L*=Lr=Ll}
(L^{*}(1)+L^{*}(0))f=L_r^{[-1]}(0)f-L_l^{[1]}(0)f,
\end{eqnarray}
which is equivalent to
\begin{eqnarray}
\<(L_r^{[-1]}(0)-L_l^{[1]}(0))f,w\>=\<f, (L(-1)+L(0))w\>
\end{eqnarray}
for $w\in W$.
Furthermore, let $f\in \mathfrak{D}_{(-1)}(W)\subset W^{*}$. 
Then $f\in (W/O'_{n,m}(W))^{*}$ with $m,n\in \N$ if and only if $f\in \Omega_{m,n}(\mathfrak{D}_{(-1)}(W))$ and
\begin{eqnarray}
 L_r^{[-1]}(0)f-L_l^{[1]}(0)f=(n-m)f.
\end{eqnarray}
\ep

\begin{proof} Recall from Proposition \ref{tri-relation} that for $v\in V,\ f\in  \mathfrak{D}_{(-1)}(W)$, we have
\begin{eqnarray*}
x_0^{-1}\delta\left(\frac{x+1}{x_0}\right)Y_W^{*}(v,x)f-x_0^{-1}\delta\left(\frac{1+x}{x_0}\right)Y^{R}(v,x)f
=x^{-1}\delta\left(\frac{-1+x_0}{x}\right)Y^{L}(v,x_0)f.\  \
\end{eqnarray*}
Specializing $v$ to the conformal vector $\omega$ and applying $\Res_{x_0}(x+x^2)$, we get
\begin{eqnarray*}
(x+x^2)Y_W^{*}(\omega,x)f- (x+x^2)Y^R(\omega,x)f=\Res_{x_0}x^{-1}\delta\left(\frac{-1+x_0}{x}\right)(x+x^2)Y^L(\omega,x_0)f,
\end{eqnarray*}
which by applying $\Res_x$ to both sides yields 
\begin{eqnarray}
(L^{*}(0)+L^{*}(1))f- \Res_x (x+x^2)Y^R(\omega,x)f =-\Res_{x_0}x_0(1-x_0)Y^L(\omega,x_0)f.
\end{eqnarray}
As $L(1)\omega=0$, from (\ref{degree-0-formula}) we have
\begin{eqnarray}
&&L_r^{[-1]}(0)f=\Res_x x(1+x)Y^R(\omega,x)f,\\
&&L_l^{[1]}(0)f=\Res_x x(1-x)Y^L(\omega,x)f.
\end{eqnarray}
Then we immediately get (\ref{L*=Lr=Ll}).

Let $g\in W^{*}$. Then $g=0$ on $(L(-1)+L(0)+m-n)W$ if and only if 
\begin{eqnarray}
(L^{*}(1)+L^{*}(0))g=(n-m)g\ \text{ in } W^*.
\end{eqnarray}
Furthermore, by Proposition \ref{p-A-D-space} and the first assertion 
we see that $g=0$ on $O'_{n,m}(W)$  if and only if 
$g\in\Omega_{m,n}( \mathfrak{D}_{(-1)}(W))$ and $L_r^{[-1]}(0)g-L_l^{[1]}(0)g=(n-m)g$. This proves the second assertion.
\end{proof}

With $O^{\dagger}_{n,m}(W)\subset O'_{n,m}(W)$,  the identity operator on $W$ yields a linear epimorphism 
from $A^{\dagger}_{n,m}(W)$ to $A^{\diamond}_{n,m}(W)$,
and hence $(A^{\diamond}_{n,m}(W))^{*}$ is naturally a subspace of $(A^{\dagger}_{n,m}(W))^{*}$.

\bc{coro-O-m-n-prime}
Let $m,n\in \N$. Then $(A^{\diamond}_{n,m}(W))^{*}$ is a submodule of the $A_m(V)\otimes A_n(V)$-module 
$(A^{\dagger}_{n,m}(W))^{*}$ defined in Theorem \ref{thm-Amn}.
On the other hand, the image of $(L(-1)+L(0)+m-n)W$ in $A^{\dagger}_{n,m}(W)$ 
is a submodule of the $A_n(V)$-$A_m(V)$ bimodule $A^{\dagger}_{n,m}(W)$
and $A^{\diamond}_{n,m}(W)$ is naturally an $A_n(V)$-$A_m(V)$ bimodule.
\ec

\begin{proof} Recall that for any $k\in \N$, the image $[\omega]$ of the conformal vector $\omega$ in $A_k(V)$ 
is a central element and $[\omega]$ acts 
on $\Omega_k(M)$ as $L(0)$ for any weak $V$-module $M$. We see that $L_r^{[-1]}(0)-L_l^{[1]}(0)$ commutes with the action of
$A_m(V)\otimes A_n(V)$ on $\Omega_{m,n}^{[1,-1]}(\mathfrak{D}_{(-1)}(W))$.
It then follows immediately from the second part of Proposition \ref{p-L(-1)+L(0)} 
that $(A^{\diamond}_{n,m}(W))^{*}$ is an $A_m(V)\otimes A_n(V)$-submodule of $\Omega_{m,n}^{[1,-1]}(\mathfrak{D}_{(-1)}(W))$.
Then making use of the linear algebra duality again we conclude that the image of $(L(-1)+L(0)+m-n)W$ in $A^{\dagger}_{n,m}(W)$ 
is an $A_n(V)$-$A_m(V)$ submodule, so that the corresponding quotient module gives an $A_n(V)$-$A_m(V)$ bimodule structure
on $A^{\diamond}_{n,m}(W)$.
\end{proof}

\bd{def-Omega'}
{\em For $m,n\in \N$, set
\begin{eqnarray}
 \Omega^{\diamond}_{m,n}(\mathfrak{D}_{(-1)}(W))=\{ f\in \Omega_{m,n}(\mathfrak{D}_{(-1)}(W))\ |\ (L_r^{[-1]}(0)-L_l^{[1]}(0))f=(n-m)f\}.
 \end{eqnarray}}
 \ed
 
 By Proposition \ref{p-L(-1)+L(0)}, we have
 \begin{eqnarray}\label{O-diamond-A-diamond}
 \Omega^{\diamond}_{m,n}(\mathfrak{D}_{(-1)}(W))=(A^{\diamond}_{n,m}(W))^{*}.
 \end{eqnarray}

From Definition \ref{def-Omega'}, for any fixed $m\in \N$, the sum $\sum_{n\in \N}\Omega^{\diamond}_{m,n}(\mathfrak{D}_{(-1)}(W))$ 
in $\mathfrak{D}_{(-1)}(W)$ is a direct sum. For $m\in \N$, set
\begin{eqnarray}
\Omega^{\diamond}_{m,\Box}(\mathfrak{D}_{(-1)}(W)) =\bigoplus_{n\in \N}\Omega^{\diamond}_{m,n}(\mathfrak{D}_{(-1)}(W))
 \subset \mathfrak{D}_{(-1)}(W).
\end{eqnarray}

Note that for homogeneous vector $u\in V$ and for $p\in \Z$,
$$[L_l^{[1]}(0), u^{R[-1]}_{\wt u-1+p}]=0\ \ \text{and }\ [L_r^{[-1]}(0), u^{R[-1]}_{\wt u-1+p}]=-pu^{R[-1]}_{\wt u-1+p}.$$
 Combining this with Lemma \ref{lemma-omega-union}, we immediately have:

\bp{lem-O-diamond-mBox}
The subspace $\Omega^{\diamond}_{m,\Box}(\mathfrak{D}_{(-1)}(W))$ of $\mathfrak{D}_{(-1)}(W)$
 is an $\N$-graded weak (sub)module of $(\mathfrak{D}_{(-1)}(W), (Y^R)^{[-1]})$ (and $(\mathfrak{D}_{(-1)}(W), Y^{R})$).
Furthermore, we have
\begin{eqnarray}
u^{R[-1]}_{\wt u-1+p}\cdot \Omega^{\diamond}_{m,n}(\mathfrak{D}_{(-1)}(W))\subset \Omega^{\diamond}_{m,n-p}(\mathfrak{D}_{(-1)}(W))
\end{eqnarray}
for any homogeneous vector $u\in V$ and for any $p,n\in \Z$.
\ep

Next, we consider the dual case. Using (\ref{other-direction}), (\ref{O-diamond-A-diamond}), and Proposition \ref{lem-O-diamond-mBox}, 
by a straightforward argument we have the following result of \cite{dj1} (cf. Corollary \ref{coro-O-m-n}): 

\bl{v[p]-O'}
Let $v\in V,\ p\in \Z,\ m,n\in \N$. Then
\begin{eqnarray}
v[p]\bar{*}_{m}^{n}O'_{n,m}(W)\subset O'_{n+p,m}(W),
\end{eqnarray}
where by definition $O'_{k,m}(W)=W$ for $k<0$ (so that $A_{k,m}^{\diamond}(W)^{*}=0=\Omega_{m,k}^{\diamond}(\mathfrak{D}_{(-1)}(W))$).
\el

The following definition was due to \cite{dj1}:

\bd{def-A-diamond}
{\em Let $m\in \N$. Form an $\N$-graded vector space
\begin{eqnarray}
A^{\diamond}_{\Box,m}(W)=\bigoplus_{n\in \N}A^{\diamond}_{n,m}(W),
\end{eqnarray}
which is naturally a right $A_m(V)$-module.}
\ed

From (\ref{O-diamond-A-diamond}), $\Omega^{\diamond}_{m,\Box}(\mathfrak{D}_{(-1)}(W))$ is the graded dual space of $A^{\diamond}_{\Box,m}(W)$.
In the following we equip $A^{\diamond}_{\Box,m}(W)$ with an $\N$-graded weak $V$-module as the contragredient module.

\bd{def-v[p]-A-diamond}
{\em Let $m\in \N$. For $v\in V,\ p\in \Z$, 
define a homogeneous operator $v[p]$ of degree $p$  on $A^{\diamond}_{\Box,m}(W)$ by
\begin{eqnarray}
v[p]\cdot (w+O'_{n,m}(W))=v[p]\bar{*}_{m}^{n}w +O'_{n+p,m}(W)\in A^{\diamond}_{n+p,m}(W)
\end{eqnarray}
for $n\in \Z,\ w\in W$. (Recall Lemma \ref{v[p]-O'}).}
\ed

Then we have the following result of \cite{dj1} (Remark 4.15) with a different proof:

\bt{thm-last}
Let $W$ be a weak $V$-module and let $m\in \N$. Define a linear map
$$Y^{\diamond}(\cdot,x):\ V\rightarrow (\End A^{\diamond}_{\Box,m}(W))[[x,x^{-1}]]$$
by the condition that for homogeneous vector $v\in V$, 
\begin{eqnarray}
Y^{\diamond}(v,x)=\sum_{p\in \Z}v[p]x^{p-\wt v}.
\end{eqnarray}
Then $(A^{\diamond}_{\Box,m}(W), Y^{\diamond})$ is a weak $V$-module, which coincides with the contragredient module
of $\Omega^{\diamond}_{m,\Box}(\mathfrak{D}_{(-1)}(W))$. Furthermore, the action of $V$ commutes with the action of $A_m(V)$.
\et

\begin{proof} Let $f\in \Omega^{\diamond}_{m,n}(\mathfrak{D}_{(-1)}(W))=(A^{\diamond}_{n,m}(W))^{*}$ with $m,n\in \N$ 
and let $v\in V$. We claim
\begin{eqnarray}\label{claim}
\< (Y^R)^{[-1]}(v,x)f,w\>=\<f, Y^{\diamond}(e^{xL(1)}(-x^{-2})^{L(0)}v,x^{-1})w\>
\end{eqnarray}
for $w\in W$. Assume that $v$ is homogeneous. Using Lemma \ref{lem-general-YR} we get
\begin{eqnarray*}
&&\<f,Y^{\diamond}(e^{zL(1)}(-z^{-2})^{L(0)}v,z^{-1})w\>\\
&=&\sum_{r\ge 0}\frac{1}{r!}z^{r-2\wt v}\<f,Y^{\diamond}(L(1)^r(-1)^{L(0)}v,z^{-1})w\>\\
&=&\sum_{r\ge 0}\frac{1}{r!}z^{r-2\wt v}\sum_{p\in \Z}\<f,(L(1)^r(-1)^{L(0)}v)[p]\bar{*}_{m}^nw\> z^{\wt (L(1)^rv)-p} \nonumber\\   
&=&\sum_{r\ge 0}\sum_{p\in \Z}\sum_{i=0}^{n+|p|}\frac{1}{r!}z^{r-2\wt v}\binom{-m-p-1}{i}\nonumber\\
&&\hspace{0.5cm}  
\times \Res_x \frac{(1+x)^{\wt (L(1)^rv)+m}} {x^{p+m+i+1}} \<f, Y_W(L(1)^r(-1)^{L(0)}v,x)w\> z^{\wt (L(1)^rv)-p}\\
&=&\sum_{p\in \Z}\<f, \theta(v)[p]\bar{*}_{m}^nw\>z^{-\wt v-p} \\
&=&\sum_{p\in \Z}\<v^{R[-1]}_{\wt v-1+p}f, w\>z^{-\wt v-p} \\
&=&\< (Y^R)^{[-1]}(v,z)f,w\>.
\end{eqnarray*}
This proves (\ref{claim}). Then it follows from the result of \cite{fhl} formulated as Lemma \ref{lemma-graded-dual} below.
 \end{proof}

The following is a variation of Theorem 5.2.1 and Proposition 5.3.1 \cite{fhl}: 

\bl{lemma-graded-dual}
Let $V$ be a vertex operator algebra and let $E=\bigoplus_{n\in \N}E(n)$ be an $\N$-graded vector space equipped
 with a linear map 
 $$Y_E(\cdot,x): V\rightarrow (\End E)[[x,x^{-1}]];\quad v\mapsto Y_E(v,x)=\sum_{p\in \Z}v_px^{-p-1}.$$ 
 Let $E'=\bigoplus_{n\in \N}E(n)^{*}$ be the graded dual space of $W$ and let 
 $$Y_E^{*}(\cdot,x): \ V\rightarrow (\End E')[[x,x^{-1}]];\ \  v\mapsto Y_E^{*}(v,x)$$ 
 be a linear map such that
 $$\<Y_E^{*}(v,x)f,w\>=\<f, Y_E(e^{xL(1)}(-x^{-2})^{L(0)}v,x^{-1})w\>$$
 for $v\in V,\ f\in E',\ w\in E$.
  Then $(E', Y_E^*)$ is an $\N$-graded weak $V$-module if and only if 
 $(E,Y_E)$ is an $\N$-graded weak $V$-module.
 \el

In the following, we discuss one of the main theorems in \cite{dj1} in the revised setup.
Let $m\in \N$ and let $U$ be an $A_m(V)$-module. Follow \cite{dj1} to define
\begin{eqnarray}
A^{\diamond}_{\Box,m}(V)\otimes _{A_m(V)}U=\bigoplus_{n\in \N}A^{\diamond}_{n,m}(V)\otimes _{A_m(V)}U,
\end{eqnarray}
which is an $\N$-graded weak $V$-module. The following is part of Theorem 4.13 in \cite{dj1} with a suitably modified proof:

\bl{lemma-U}
We have $A^{\diamond}_{m,m}(V)\otimes_{A_m(V)}U=U$
as (left) $A_m(V)$-modules. Furthermore, $U$ generates $A^{\diamond}_{\Box,m}(V)\otimes _{A_m(V)}U$ as a weak $V$-module.
\el

\begin{proof} From definition, we see that $O'_{m,m}(V)=O_m(V)$ and hence $A^{\diamond}_{m,m}(V)=A_m(V)$ as vector spaces.
Let $u,v\in V$ be homogeneous vectors. By Lemma \ref{left-right-diff} we have
 \begin{eqnarray*}
 u\bar{*}_{m}^nv-u\bar{*}_{m,n}v=\Res_x (1+x)^{\wt u-1}Y(u,x)v.
 \end{eqnarray*}  
 On the other hand, recall from \cite{dlm-anv} that  
 $$ u*_{m}v-v*_{m}u-\Res_x (1+x)^{\wt u-1}Y(u,x)v\in O_m(V).$$
 From definition, $u\bar{*}_{m}^mv=u*_mv$, so $u\bar{*}_{m,n}v-v*_mu\in O_m(V)$.
 Thus $A^{\diamond}_{m,m}(V)=A_m(V)$ as $A_m(V)$-bimodules.
Therefore,  we have
$$A^{\diamond}_{m,m}(V)\otimes_{A_m(V)}U=A_{m}(V)\otimes_{A_m(V)}U=U$$
as (left) $A_m(V)$-modules. 

Let $n\in \N$ and let $u\in V$ be any homogeneous vector. From (\ref{deform-dj}) we have
\begin{eqnarray*}
u[n-m]\bar{*}_{m}^{m}{\bf 1}&=&\Res_x \sum_{i=0}^{m+|n-m|}\binom{-n-1}{i}
\frac{(1+x)^{m+\wt u}} {x^{n+i+1}}Y(u,x){\bf 1}\\
&=&\Res_x \sum_{i=0}^{m+|n-m|}\binom{-n-1}{i}
\frac{(1+x)^{m+\wt u}} {x^{n+i+1}}e^{xL(-1)}u\\
&=&\Res_x \sum_{i=0}^{m+|n-m|}\sum_{j,k\ge 0}\binom{-n-1}{i}
\binom{m+\wt u}{j} \frac{x^{j+k}}{x^{n+i+1}}\frac{1}{k!}L(-1)^{k}u.
\end{eqnarray*}
Noticing that for any homogeneous vector $v\in V$, 
$$L(-1)v\equiv (n-m-\wt v)v\ \ \mod\; O'_{n,m}(V),$$
we have
\begin{eqnarray*}
\frac{1}{k!}L(-1)^{k}u&\equiv &\frac{1}{k!}(n-m-(\wt u+k-1))\cdots (n-m-\wt u)u\ \ \mod\; O'_{n,m}(V)\\
&=&\binom{n-m-\wt u}{k}u
\end{eqnarray*}
for $k\ge 0$. Then
\begin{eqnarray*}
u[n-m]\bar{*}_{m}^{m}{\bf 1}\equiv \sum_{i=0}^{m+|n-m|}\sum_{j,k\ge 0}\binom{-n-1}{i}\binom{m+\wt u}{j}\binom{n-m-\wt u}{k}u
\ \ \mod\; O'_{n,m}(V),
\end{eqnarray*}
where $n+i=j+k$ is assumed. Note that for $r\ge n$, we have
$$\sum_{j,k\ge 0,\; j+k=r}\binom{m+\wt u}{j}\binom{n-m-\wt u}{k}=\binom{n}{r}=\delta_{r,n}.$$
Thus $u[n-m]\bar{*}_{m}^{m}{\bf 1}\equiv u\ \ \mod\; O'_{n,m}(V)$. That is, $u[n-m]({\bf 1}+O'_{m,m}(V))=u+O'_{n,m}(V)$.
This shows that $U$ generates $A^{\diamond}_{\Box,m}(V)\otimes_{A_m(V)}U$ as a weak $V$-module.
\end{proof}

The following is Theorem 4.13 in \cite{dj1} with a suitably modified proof:

\bt{A-nm-W}
Let $m\in \N$ and let $U$ be an $A_m(V)$-module.  
Suppose that $W$ is a weak $V$-module and $\psi: U\rightarrow \Omega_m(W)$ is an $A_m(V)$-module morphism. 
For $n\in \N$, define a bilinear map $F_{n,m}:\ V\times U\rightarrow W$ by
\begin{eqnarray}
F_{n,m}(v,w)=\Res_x x^{m-n-1} Y(x^{L(0)}v,x)\psi(w)
\end{eqnarray}
for $v\in V,\ w\in U$. Then $F_{n,m}(v,w)\in \Omega_n(W)$,
$F_{n,m}(O'_{n,m}(V)\times U)=0$, and
$F_{n,m}$ reduces to an $A_n(V)$-module morphism 
$F_{n,m}^{\diamond}:\ A^{\diamond}_{n,m}(V)\otimes _{A_m(V)}U\rightarrow \Omega_n(W)$. Define a linear map
$$\tilde{\psi}:\ A^{\diamond}_{\Box,m}(V)\otimes _{A_m(V)}U\rightarrow W$$
by $\tilde{\psi}|_{A^{\diamond}_{n,m}(V)\otimes _{A_m(V)}U}=F^{\diamond}_{n,m}$ for $n\in \N$.
Then $\tilde{\psi}$ is a $V$-module morphism which is uniquely determined by the condition $\tilde{\psi}|_U=\psi$.
\et

\begin{proof}  For $v\in V,\ w\in U$, as $\psi(w)\in \Omega_m(W)$, by Lemma \ref{lemma-omega-union}, 
we have $F_{n,m}(v,w)\in \Omega_n(W)$. 
Let $u,v\in V$ be homogeneous vectors and let $w\in U$. Recall 
$$u\circ_{m}^nv=\Res_z \frac{(1+z)^{\wt u+m}}{z^{m+n+2}}Y(u,z)v.$$
Noticing that $x^{\wt u+m}Y(u,x)\psi(w)\in W[[x]]$, 
from Jacobi identity, we have 
\begin{eqnarray}\label{weak-assoc-Y0}
(x_0+x_2)^{\wt u+m}Y(u,x_0+x_2)Y(v,x_2)\psi(w)=(x_2+x_0)^{\wt v+m}Y(Y(u,x_0)v,x_2)\psi(w).
\end{eqnarray}
Using this relation we get
\begin{eqnarray*}
&&F_{n,m}(u\circ_{m}^nv,w)\\
&=&\Res_x\Res_z x^{m-n-1}\frac{(1+z)^{\wt u+m}}{z^{m+n+2}} Y\left(x^{L(0)}Y(u,z)v,x\right)\psi(w)\\
&=&\Res_x\Res_z x^{m-n-1+\wt u+\wt v}\frac{(1+z)^{\wt u+m}}{z^{m+n+2}} Y\left(Y(u,xz)v,x\right)\psi(w)\\
&=&\Res_{z_0}\Res_{x}x^{\wt v+m}\frac{(x+z_0)^{\wt v+m}}{z_0^{m+n+2}}Y(Y(u,z_0)v,x)\psi(w)\\
&=&\Res_{z_0}\Res_{x}x^{\wt v+m}\frac{(z_0+x)^{\wt v+m}}{z_0^{m+n+2}}Y(u,z_0+x)Y(v,x)\psi(w)\\
&=&0
\end{eqnarray*}
as $x^{\wt v+m}Y(v,x)\psi(w)\in W[[x]]$. On the other hand, we have
\begin{eqnarray*}
&&F_{n,m}((L(-1)+L(0)+m-n)v,w)\\
&=&\Res_x x^{m-n-1}Y\left(x^{L(0)}(L(-1)+L(0)+m-n)v,x\right)\psi(w)\\
&=&\Res_x x^{m-n+\wt v}Y(L(-1)v,x)\psi(w)
+x^{m-n-1+\wt v}Y((L(0)+m-n)v,x)\psi(w)\\
&=&\Res_x x^{m-n+\wt v}\frac{d}{dx}Y(v,x)\psi(w)+\Res_x (\wt v+m-n)x^{m-n-1+\wt v}Y(v,x)\psi(w)\\
&=&\Res_x \frac{d}{dx}\left(x^{m-n+\wt v}Y(v,x)\psi(w)\right)\\
&=&0.
\end{eqnarray*}
This shows that $F_{n,m}$ vanishes on $O'_{n,m}(V)\times U$.

Noticing that as $\psi(w)\in \Omega_m(W)$,
$$x^{m-n-1+\wt v+i}Y(v,x)\psi(w)\in W[[x]]$$
for $i>n$, using Lemma \ref{left-right-diff} we have
\begin{eqnarray*}
&&F_{n,m}(u\bar{*}_{m,n}v,w)\\
&=&F_{n,m}(u\bar{*}_{m}^{n}v,w)-\Res_z(1+z)^{\wt u-1}F_{n,m}(Y(u,z)v,w)\\
&=&\Res_{x}\Res_{z}\sum_{i=0}^{n}\binom{-m-1}{i}\frac{(1+z)^{\wt u+m}}{z^{m+i+1}}x^{m-n-1}Y(x^{L(0)}Y(u,z)v,x)\psi(w)\\
&&-\Res_{x}\Res_{z}(1+z)^{\wt u-1}x^{m-n-1}Y(x^{L(0)}Y(u,z)v,x)\psi(w)\\
&=&\Res_{x}\Res_{z}\sum_{i=0}^{n}\binom{-m-1}{i}\frac{(1+z)^{\wt u+m}}{z^{m+i+1}}x^{m-n-1+\wt u+\wt v}Y(Y(u,xz)v,x)\psi(w)\\
&&-\Res_{x}\Res_{z}x^{m-n-1+\wt u+\wt v}(1+z)^{\wt u-1}Y(Y(u,xz)v,x)\psi(w)
\end{eqnarray*}
and
\begin{eqnarray*}
&&F_{n,m}(v,u\cdot w)\\
&=&\Res_{x_1}\Res_{x}x_1^{\wt u-1} x^{m-n-1+\wt v}Y(v,x)Y(u,x_1)\psi(w)\\
&=&\Res_{x_1}\Res_{x}x_1^{\wt u-1} x^{m-n-1+\wt v}Y(u,x_1)Y(v,x)\psi(w)\\
&&-\Res_{x_0}\Res_{x} x^{m-n-1+\wt v}(x+x_0)^{\wt u-1} Y(Y(u,x_0)v,x)\psi(w)\\
&=&\Res_{x_0}\Res_{x}x^{m-n-1+\wt v}(x_0+x)^{\wt u-1} Y(u,x_0+x)Y(v,x)\psi(w)\\
&&-\Res_{z}\Res_{x}x^{m-n-1+\wt u+\wt v}(1+z)^{\wt u-1} Y(Y(u,xz)v,x)\psi(w)\\
&=&\sum_{i=0}^n\binom{-m-1}{i}\Res_{x_0}\Res_{x} x_0^{-m-1-i}x^{m-n-1+\wt v+i}\\
&&\quad \cdot [(x_0+x)^{\wt u+m}Y(u,x_0+x)Y(v,x)\psi(w)]\\
&&-\Res_{z}\Res_{x}x^{m-n-1+\wt u+\wt v}(1+z)^{\wt u-1} Y(Y(u,xz)v,x)\psi(w)\\
&=&\sum_{i=0}^n\binom{-m-1}{i}\Res_{x_0}\Res_{x} x_0^{-m-1-i}x^{m-n-1+\wt v+i}(x+x_0)^{\wt u+m}Y(Y(u,x_0)v,x)\psi(w)\\
&&-\Res_{z}\Res_{x}x^{m-n-1+\wt u+\wt v}(1+z)^{\wt u-1} Y(Y(u,xz)v,x)\psi(w)\\
&=&\sum_{i=0}^n\binom{-m-1}{i}\Res_{z}\Res_{x} x^{m-n-1+\wt u+\wt v}\frac{(1+z)^{\wt u+m}}{z^{m+i+1}}Y(Y(u,xz)v,x)\psi(w)\\
&&-\Res_{z}\Res_{x}x^{m-n-1+\wt u+\wt v}(1+z)^{\wt u-1} Y(Y(u,xz)v,x)\psi(w).
\end{eqnarray*}
Consequently, we get
\begin{eqnarray}\label{left-right-comm}
F_{n,m}(u\bar{*}_{m,n}v,w)=F_{n,m}(v,u\cdot w).
\end{eqnarray}
On the other hand, we have
\begin{eqnarray}\label{4terms-relation}
&&F_{n,m}(v,u\cdot w)\nonumber\\
&=&\Res_{x_1}\Res_{x}x_1^{\wt u-1} x^{m-n-1+\wt v}Y(v,x)Y(u,x_1)\psi(w)\nonumber\\
&=&\Res_{x_1}\Res_{x}x_1^{\wt u-1} x^{m-n-1+\wt v}Y(u,x_1)Y(v,x)\psi(w)\nonumber\\
&&-\Res_{x_0}\Res_{x} x^{m-n-1+\wt v}(x+x_0)^{\wt u-1}Y(Y(u,x_0)v,x)\psi(w)\nonumber\\
&=&u\cdot F_{n,m}(v,w) -\Res_{z}\Res_{x} x^{m-n-1+\wt u+\wt v}(1+z)^{\wt u-1}Y(Y(u,xz)v,x)\psi(w)\nonumber\\
&=&u\cdot F_{n,m}(v,w)
 -\Res_{z}\Res_{x} x^{m-n-1}(1+z)^{\wt u-1} Y(x^{L(0)}Y(u,z)v,x)\psi(w)\nonumber\\
&=&u\cdot F_{n,m}(v,w)+F_{n,m}(u\bar{*}_{m,n}v,w)-F_{n,m}(u\bar{*}_{m}^nv,w).
\end{eqnarray}
Combining this with (\ref{left-right-comm}) we get
\begin{eqnarray}
u\cdot F_{n,m}(v,w)=F_{n,m}(u\bar{*}_{m}^nv,w).
\end{eqnarray}
Therefore, $F_{n,m}$ reduces to an $A_n(V)$-module morphism $F_{n,m}^{\diamond}:\ A^{\diamond}_{n,m}(V)\otimes _{A_m(V)}U\rightarrow W$. 

Now, using $F^{\diamond}_{n,m}$ for all $n\in \N$ we get a linear map 
$$\tilde{\psi}:\  A^{\diamond}_{\Box,m}(V)\otimes _{A_m(V)}U\rightarrow W.$$
Let $u,v\in V$ be homogeneous and let $w\in U,\ p\in \Z$. Noticing that
$$x^{m-n-1+i}Y(x^{L(0)}v,x)\psi(w)\in W[[x]]$$
for all $i>n+|p|\ (\ge n)$, we have
\begin{eqnarray*}
&&\Res_{x_1}x_1^{\wt u-1-p}Y(u,x_1)F_{n,m}(v,w)\\
&=&\Res_{x_1}\Res_{x}x_1^{\wt u-1-p}x^{m-n-1}Y(u,x_1)Y(x^{L(0)}v,x)\psi(w)\\
&=&\Res_{x_0}\Res_{x}(x_0+x)^{\wt u-1-p}x^{m-n-1}Y(u,x_0+x)Y(x^{L(0)}v,x)\psi(w)\\
&=&\Res_{x_0}\Res_{x}\sum_{i=0}^{n+|p|}\binom{-m-p-1}{i}x_0^{-m-p-1-i}x^{m-n-1+i}\\
&&\quad \quad \cdot (x_0+x)^{\wt u+m}Y(u,x_0+x)Y(x^{L(0)}v,x)\psi(w)\\
&=&\Res_{x_0}\Res_{x}\sum_{i=0}^{n+|p|}\binom{-m-p-1}{i}x_0^{-m-p-1-i}x^{m-n-1+i}\\
&&\quad \quad \cdot (x+x_0)^{\wt u+m}Y(Y(u,x_0)x^{L(0)}v,x)\psi(w)\\
&=&\Res_{x_0}\Res_{x}\sum_{i=0}^{n+|p|}\binom{-m-p-1}{i}x_0^{-m-p-1-i}x^{m-n-1+i}\\
&&\quad \quad \cdot (x+x_0)^{\wt u+m}Y(x^{L(0)}Y(x^{-L(0)}u,x_0/x)v,x)\psi(w)\\
&=&\Res_{z}\Res_{x}\sum_{i=0}^{n+|p|}\binom{-m-p-1}{i}x^{m-n-1-p}\frac{(1+z)^{\wt u+m}}{z^{m+p+1+i}}Y(x^{L(0)}Y(u,z)v,x)\psi(w)\\
&=&\Res_{x}x^{m-n-1-p}Y(x^{L(0)}(u[p]\bar{*}_m^nv),x)\psi(w)\\
&=&\Res_{x_1}x_1^{\wt u-1-p}F_{n+p,m}\left(Y^{\diamond}(u,x_1)(v+O'_{n,m}(V)),w\right).
\end{eqnarray*}
This proves that $\tilde{\psi}$ is a $V$-module homomorphism. The uniqueness is clear.
\end{proof}

\end{document}